\documentclass{article}[12pt]
\usepackage{mathrsfs}
\usepackage{mathtools}
\usepackage{graphicx}
\usepackage{epstopdf}
\usepackage{stmaryrd}
 \usepackage{subfigure}
\usepackage{float}
\usepackage{amsmath}
    \usepackage{bm}
\usepackage{amsfonts,amssymb}
  \usepackage{dsfont}
\usepackage{amsfonts}
\usepackage{mathrsfs}
  \usepackage{pifont}
\numberwithin{equation}{section}
 \usepackage{amsthm}

 \usepackage{amsmath}
 \usepackage{hyperref}
 \usepackage{multirow}

\begin{document} 
\title{{\large  \textbf{Decoupling Numerical Method Based on Deep Neural Network for Nonlinear Degenerate Interface Problems}}}

\author{Chen Fan$^{a}$, Zhiyue Zhang$^{a,}$\footnote{Corresponding author. E-mail address: zhangzhiyue@njnu.edu.cn.}\\
{\small ${a}$ School of Mathematical Sciences, Jiangsu Key Laboratory for NSLSCS,} \\
{\small Nanjing Normal University, Nanjing 210023, China }} 

\date{}
\maketitle 
{\footnotesize
\noindent \textbf{Abstract~~~} Interface problems depict many fundamental physical phenomena and widely apply in the engineering. However, it is challenging to develop efficient fully decoupled numerical methods for solving degenerate interface problems in which
the coefficient of a PDE is discontinuous and greater than or equal to zero on the interface. The main motivation in this paper is to construct fully decoupled numerical methods for solving nonlinear degenerate interface problems with ``double singularities". An efficient fully decoupled numerical method is proposed for nonlinear degenerate interface problems. The scheme combines deep neural network on the singular subdomain with finite difference method on the regular subdomain. The key of the new approach is to split nonlinear degenerate partial differential equation with interface into two independent boundary value problems based on deep learning. The outstanding advantages of the proposed schemes are that not only the convergence order of the degenerate interface problems on whole domain is determined by the finite difference scheme on the regular subdomain, but also can calculate $\mathbf{VERY}$ $\mathbf{BIG}$ jump ratio(such as $10^{12}:1$ or $1:10^{12}$) for the interface problems including degenerate and non-degenerate cases. The expansion of the solutions does not contains any undetermined parameters in the numerical method. In this way, two independent nonlinear systems are constructed in other subdomains and can be computed in parallel. The flexibility, accuracy and efficiency of the methods are validated from various experiments in both 1D and 2D. Specially, not only our method is suitable for solving degenerate interface case, but also for non-degenerate interface case. Some application examples with complicated multi-connected and sharp edge interface examples including degenerate and nondegenerate cases are also presented.
\\[2mm]
\textbf{Key words~~~}nonlinear degenerate interface problems; deep neural network; fully decoupled method; very big jump ratio; convergence order; sharp edge interface
\\[2mm]
\\
\textbf{Mathematics Subject Classification~~~}34B16, 35R05, 65M85, 65N06, 68T99}

\section{Introduction}
\label{s1}
{Nonlinear degenerate interface problems can depict many fundamental physical phenomena in chemical and mechanical engineering, physics and many other applications\cite{del2017numerical,zhao2010high,leveque1994immersed,sun2014adaptive,ren2000iterative,guo2020recovering}. For the standard interface problems, it has attracted great interests in numerical computations, such as finite element method\cite{guo2020recovering,arbogast1996nonlinear,cai2017discontinuous,Chen1998Finite,xia2014mib}, finite difference method\cite{adams2002immersed, zhou2006high, beale2019solution, beale2007accuracy}, finite volume element method \cite{wang2021bilinear,ewing1999immersed, cao2017superconvergence,zhu2015immersed,wang2021new}, spectral method\cite{shen2011spectral,jiang2012phase,Chen2018Enriched}, least-squares method \cite{xu2021fourth} and references therein. There has been a great deal of rigorous mathematical theory and numerical analysis to deal with degenerate PDE\cite{wang2013approximate,wang2014carleman,du2012analysis,bao2017numerical,shen2016,wu2019finite,ben2001jacobi,bernis1990higher,gunzburger2018stokes,zhu2019fast}. To the best of our knowledge, degenerate interface problems have received less attention so far, a few notable approaches can be found in the literature to handle the degenerate PDE with interface\cite{zhao2021semi,zhou2006high,bedrossian2010second}. As is well known, the difficulty lies in the ``double singularities" for nonlinear degenerate interface problems, namely, degeneracy and interface. Generally speaking, the most expensive part work of numerical schemes on standard sharp interface problems\cite{hou2010numerical,huang2017unfitted,he2010interior,ji2022immersed} is how to approximate the jump conditions very well. For example, there are many methods are interesting but the technique to treat the jump conditions is quite complicated. Nevertheless, our proposed approach based on deep neural network uses different, simple and natural techniques to treat the singularities compared with the above references, and hence obtains numerical method to solve nonlinear degenerate interface problems. In fact, the challenge work of numerical simulation on nonlinear degenerate interface problems are how to design the numerical methods not only to reduce the singularities affect at degenerate points, but also are less dependent or independent of the jump conditions. Due to nonlinear degenerate interface problems possess ``double singularities", it is usually required extremely fine grids such as adaptive mesh or graded mesh to reduce singularity affect. Obviously, it is impossible  to use uniform grids to numerically solve nonlinear degenerate interface problems for the traditional numerical methods. The main goal of this paper is to present an efficient and fully decoupled finite difference methods with uniform grids based on deep neural network for solving nonlinear degenerate interface problems.
	
On the other hand, the deep neural network (DNN) models have achieved great successes in artificial intelligence fields including  high-dimensional problems in computer vision, natural language processing, time series analysis, pattern and speech recognition. It is worth noting that even if there is a universal approximation theoretical results about the single layer neural network, the approximation theory of DNN still remain an open question. However, this should not prevent us from trying to apply deep learning to other problems such as numerical weather forecast, petroleum engineering, turbulence flow and interface problems. There are two main techniques to solve PDEs with deep learning, the first is to parameterize the solution of PDEs by the deep neural network (DNN). One of methods is that a universal approximation based on a neural network and point collocation are used to transform the PDE into an unconstrained minimization problem. The other one is that the original problem is transformed an optimization problem with variational form based on representing the trail functions by deep neural networks. Recently, we have noticed there are some gratifying works by using mesh free methods with DNN model to solve PEDs and interface problems\cite{hu2022discontinuity,he2022mesh,wang2020mesh,baharlouei2023dnn}. However, we will use structured mesh method with deep learning to deal with degenerate interface problems  which is a challenge and is always of great interests. Although boundary conditions are absent on the singular sub-domains, which is known to be the extreme ill-posedness, it is shown that the DNN approach still has some merits in structured grids method. In addition, we use a hybrid asymptotic and augmented compact finite volume method to realize using semi-decoupling numerical method based on a uniform Cartesian mesh for solving 1D degenerate interface problem\cite{zhao2021semi}. This inspires us to develop fully decoupled numerical method for solving the degenerate PDE with interface. Although there have been a great deal of nice works for interface problems\cite{li2003new,hou2005numerical,zhao2017efficient,bedrossian2010second,zhou2006high,hou2010numerical,huang2017unfitted,he2010interior,wang2020mesh,wang2021bilinear,wang2020mesh,guo2020recovering}, there are quite a few fully decoupled numerical methods on the uniform grids for solving such interface problems, even to mentioned interesting degenerate interface problems.
	
In this paper, we focus on constructing fully decoupled numerical algorithms based on deep learning for solving the degenerate interface problems. This method not only effectively reduces the influence of the degeneracy and interface but also provides an accurate solutions on a uniform Cartesian mesh. We construct two DNN structures near the interface instead of the whole domain, and find the optimal solution by minimizing the mean squared error loss that consists of the equation and the interface conditions. These two parts are linked by its normal derivative jump conditions. We use DNN to treat considered problems on singular sub-domains near the interface to get a solution, then obtain two independent decoupled boundary value sub-problems without interface on regular sub-domains. We can compute those two nonlinear systems in parallel. We find that the proposed our approach is simple, easy to implement reducing lots efforts in handling jump conditions and also its ability to use existing method for solving nonlinear sub-problmes without interface. The choice of the singular sub-domain is more natural since we use a uniform grids, and programming of the new scheme is a straightforward task due to fully decoupled algorithms. Although deep learning has shown remarkable success in various hard problems of artificial intelligence areas, limited approximatability of deep learning with uniform grids results in two general boundary value sub-problems to get satisfactory approximations of the solutions for solving such nonlinear degenerate interface problems. A loss, no bad thing or a blessing in disguise. In fact, if deep learning has the ability to strictly decoupled the degenerate interface problems at the interface into two degenerate PDEs, we probably obtain nonlinear ill-conditioned systems for the corresponding discrete sub-problems. At this moment, we have to look for other special methods to treat degenerate PDE or interface problems likewise the litratures\cite{lusch2018deep,zhao2021semi}, and references therein.
	
The purpose of the paper is to develop a new fully decoupled numerical method based on DNN technique that not only effectively reduces the influence of the singularities and interface, but also provides a new way to realize completely decoupled method with different ideas compared to the existing methods to treat degenerate interface problems.  It does not need any extra efforts to treat the cases between degenerate interface and general interface. The proposed approach has advantages of fully decoupled two problems without interface with uniform grids. Since our fully decoupled method is independent of the interface and the jump conditions, it not only results in two independent sub-problems, but also can easily treat the cases of $\mathbf{VERY}$ $\mathbf{BIG}$ jump ratio(such as $10^{12}:1$ or $1:10^{12}$). In addition, the computational costs is almost the same for homogenous jump case and non-homogeneous jump case, this numerically demonstrates fully decoupled property of our method. The methods of this paper are sufficiently robust and also can easily handle 1D case and 2D case. In particular, it is easily to handle hard problems such as sharp-edge interface problems. Our method can robustly and efficiently apply to both of the general interface problems and degenerate interface problems, while an effective method to solve general interface problems is not suitable for solving such nonlinear degenerate interface problems. It is demonstrated that our method is a simple and straight method to deal with quit hard works. It should be mentioned that the convergence order of the schemes on entire domain for solving such degenerate PDE with interface can be determined by the convergence order of the sub-problems on regular sub-domain. Numerical experiments show that the proposed approach is able to effectively approximate the solutions of such hard degenerate interface problems. Numerical results have shown great improvement comparing to the existing methods for solving hard cases\cite{albright2017high}. From the method\cite{zhao2021semi} we know that it is impossible to split degenerate or general interface problems into two independent boundary value problems. Nevertheless, it is realized our algorithms to be completely decoupled for solving degenerate interface problems due to using dee learning. Although there are a few analytical results, the reason why deep neural networks coupled with traditional numerical methods have performed so well for solving degenerate interface problems still largely remains a mystery. This encourage us to consider the theoretical approximation analysis in the future.
	
The rest of the paper is organized as follows. In section \ref{s2}, we give some preliminaries about the Deep Neural Networks and follow this with the process on the interface and fully-decoupling two sub-problems. In section \ref{s3}, we construct Deep Neural Network structure and finite difference scheme. We present some numerical experiments including some interesting models in mathematical physics area in section \ref{s4}. Some concluding remarks are given in the final section.}


\section{Deep Neural Network}
\label{s2}
The definition and attributes of the deep neural network (DNN), particularly its approximation property, are briefly discussed in this section \cite{wang2020mesh}.

\begin{figure}
	\centering    
	\subfigure[ReLU]
	{
		\begin{minipage}{6cm}
			\centering
			\includegraphics[width=6cm,clip]{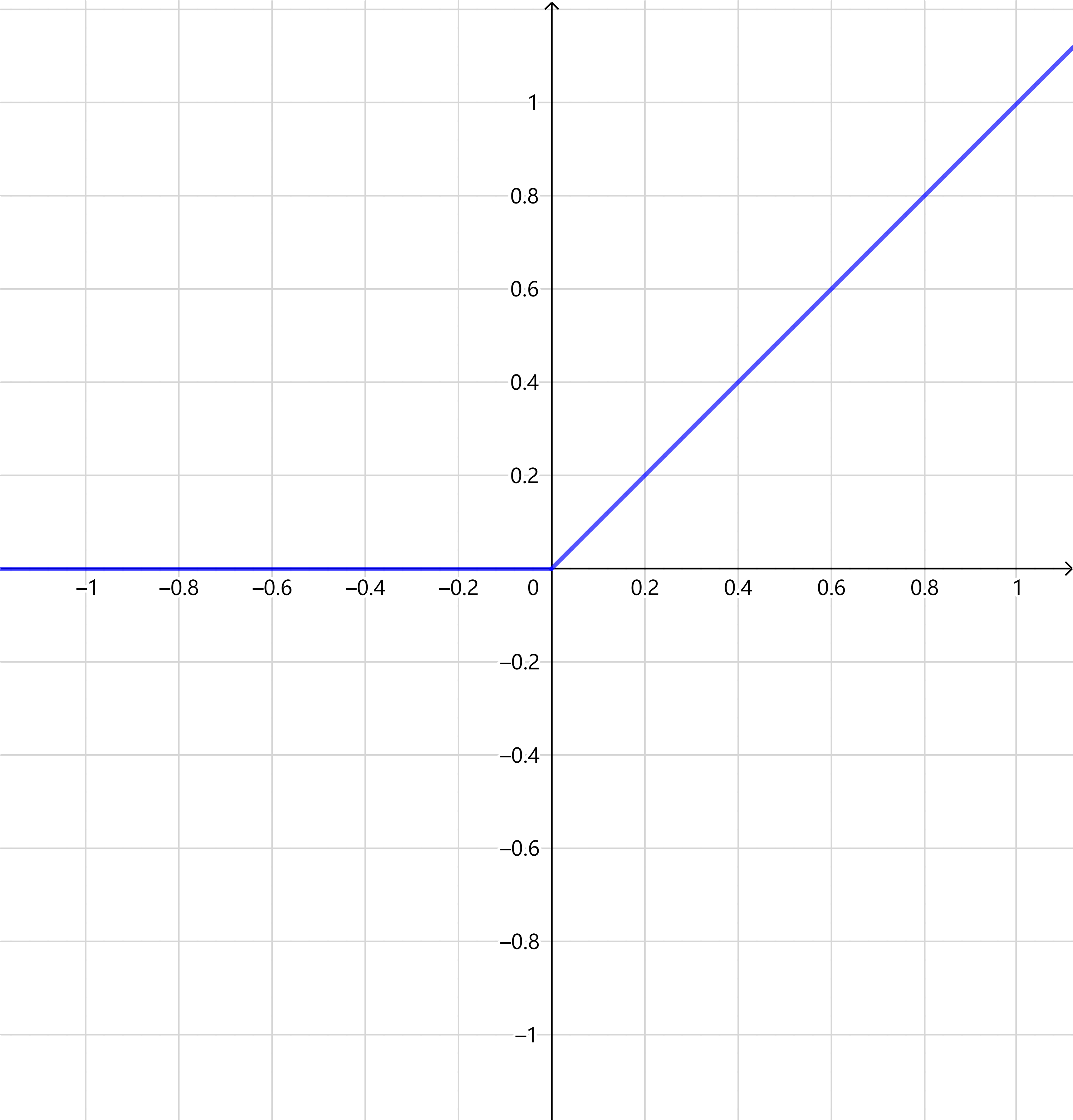}
			
		\end{minipage}
	}
	\subfigure[ELU]
	{
		\begin{minipage}{6cm}
			\centering
			\includegraphics[width=6cm,clip]{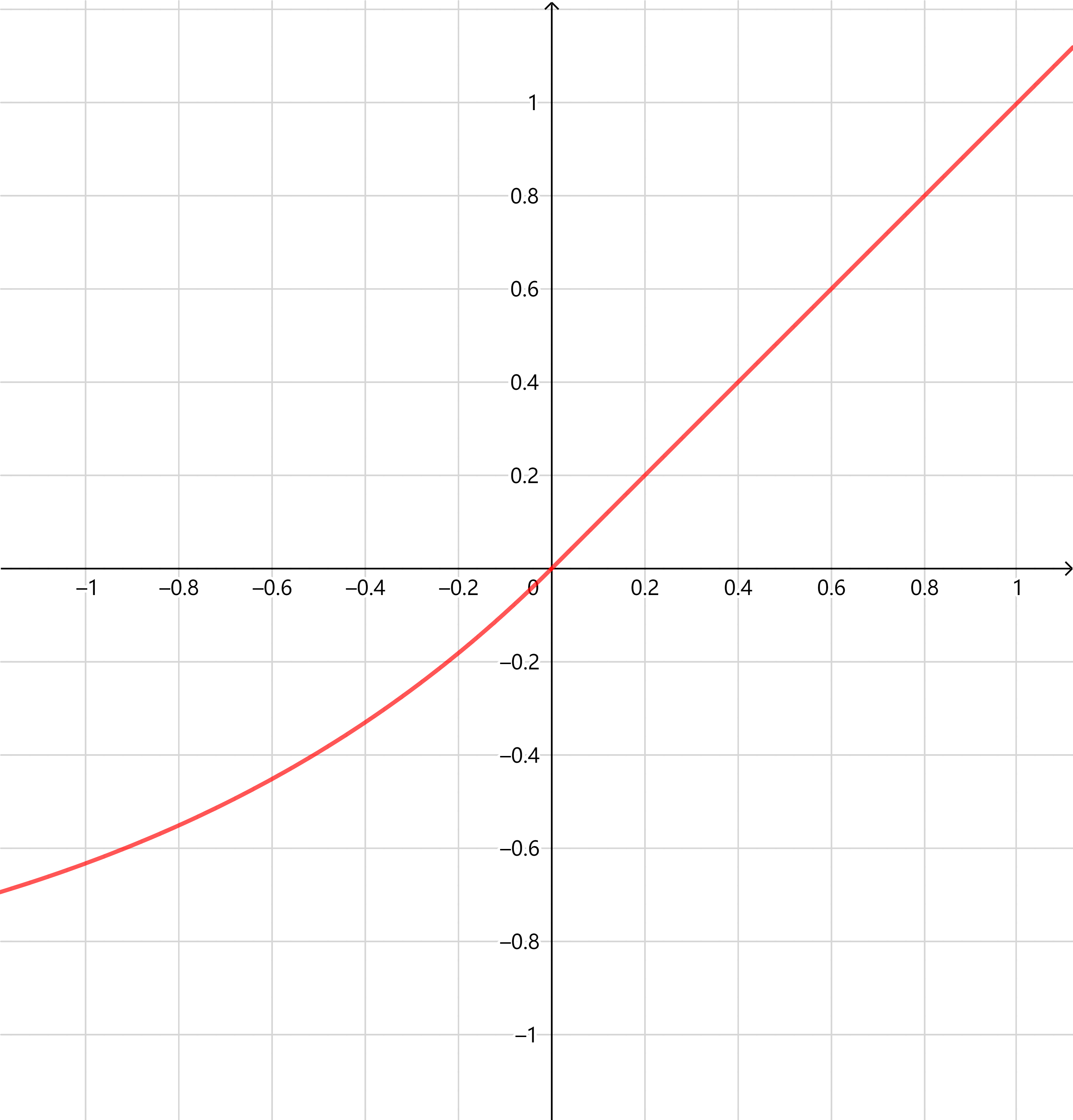}
			
		\end{minipage}
	}
	\vspace{-3mm}
	\caption{Images of the activation functions.}
	\label{figure18}
\end{figure}

In order to define a DNN, we will need two steps. The first is a (vector) linear function of the operator $T: R^n \rightarrow R^m$, defined as $T(\boldsymbol{x})=A\boldsymbol{x}+\boldsymbol{b}$, where $A=(a_{i,j}) \in R^{m \times n}$, $\boldsymbol{x}$ and $\boldsymbol{b}$ are in $R^n$ and $R^m$ respectively. A nonlinear activation function $\sigma: R \rightarrow R$ is the second. The rectified linear unit (ReLU), a commonly used activation function, is defined as $ReLU(x)=max(0,x)$\cite{lecun2015deep}. The exponential linear unit (ELU) will  be used as the activation function in this paper, defined as $ELU(x)=max(0,x)+min(0,e^x-1)$, it is mainly used to avoid the problem of gradient disappearance (Fig.\ref{figure18}). The (vector) activation function $\sigma: R^m \rightarrow R^m$ can be defined by applying the activation function in an element-wise manner.
We can define a continuous function $F(\boldsymbol{x})$ by acomposition of linear transforms and activation functions using these definitions, i.e.,
\begin{equation}
\label{e1}
F(\boldsymbol{x})=T^k \circ \sigma \circ T^{k-1} \circ \sigma \circ T^{k-2} \circ \dots \circ T^0(\boldsymbol{x}),
\end{equation}
where $T^i(\boldsymbol{x})=A_i\boldsymbol{x}+\boldsymbol{b_i}$ with $A_i$ and $\boldsymbol{b_i}$ are undetermined matrices and vectors respectively, $\sigma(x)$ being the element-wisely specified activation function to make (\ref{e1}) meaningful, the dimensions of $A_i$ and $\boldsymbol{b_i}$ were chosen. All indeterminate coefficients (e.g., $A_i$ and $\boldsymbol{b_i}$) in (\ref{e1}) are denoted as $\boldsymbol{\theta} \in \Theta $, where $\boldsymbol{\theta}$ is a high-dimensional vector and $\Theta$ is the space of $\boldsymbol{\theta}$. The DNN representation of a continuous function can be viewed as
\begin{equation}
F=F(\boldsymbol{x};\boldsymbol{\theta}).
\end{equation}
Let $\mathbb{F}=\{F(\boldsymbol{x};\boldsymbol{\theta}) \mid \boldsymbol{\theta} \in \Theta\}$ denote the set of all expressible functions by the DNN parametrized by $\boldsymbol{\theta} \in \Theta$. The approximation property of the DNN, which is relevant to the study of a DNN model's expressive power, have been discussed in other papers\cite{he2018relu,yarotsky2017error}. To accelerate the training of the neural network, we use the Adam optimizer \cite{kingma2014adam}version of the stochastic gradient descent (SGD) method in two-dimensional case\cite{robbins1951stochastic}.
 \begin{figure}
	\centering
	\includegraphics[width=15cm]{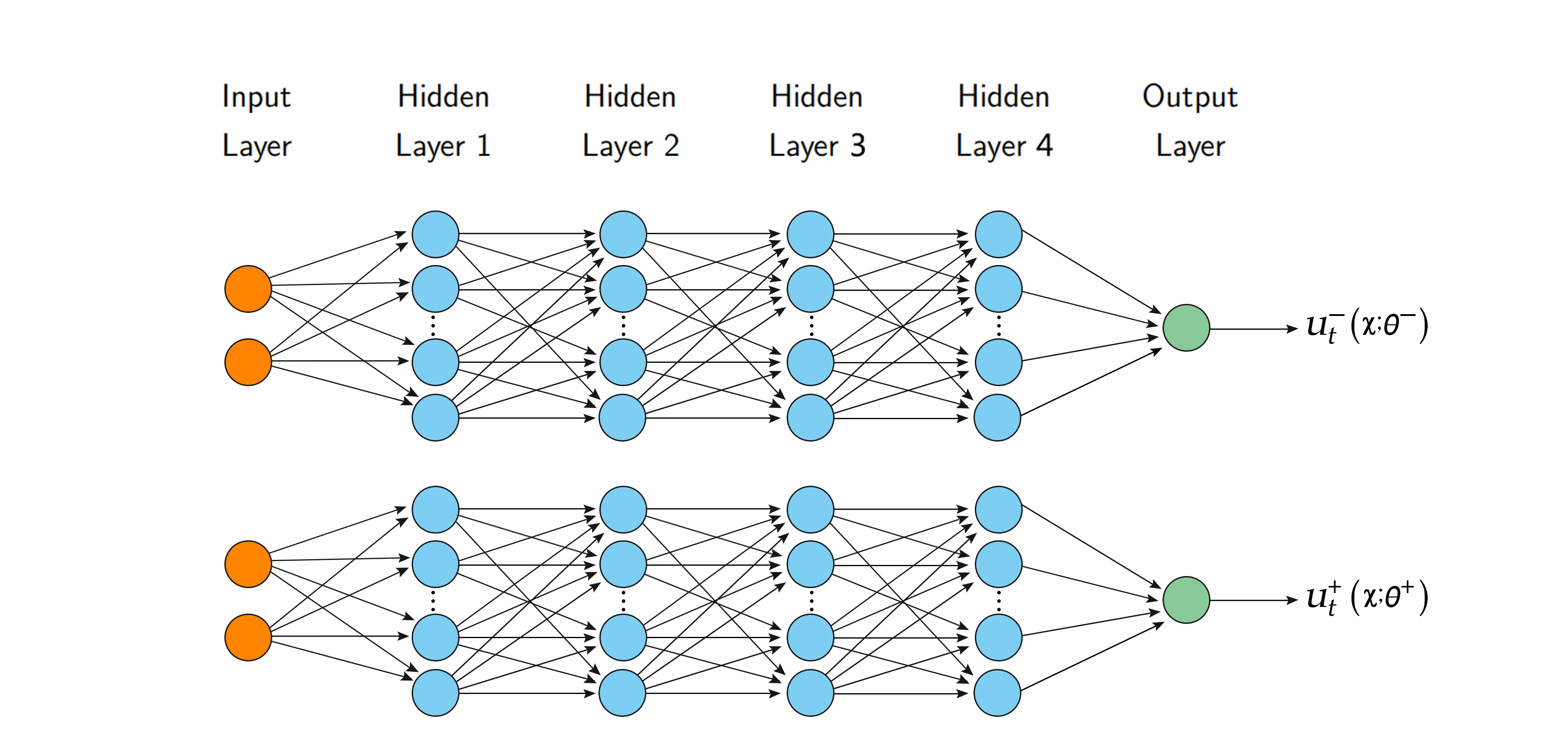}\\
	\caption{A diagram of the deep neural network architecture}
	\label{DNN}
\end{figure}


\section{2D Degenerate Elliptic Interface Problem}
\label{s3}
\subsection{Problem description}

Consider the following nonlinear degenerate elliptic equation with the interface,
\begin{equation}
\label{e2}
-\nabla \cdot(\beta(\boldsymbol{x}) \nabla u)=f(\boldsymbol{x},u), \text { in } \Omega^- \cup \Omega^+,
\end{equation}
\begin{equation*}
[u]=w, \text { on } \Gamma,
\end{equation*}
\begin{equation*}
[\beta(\boldsymbol{x}) \nabla u \cdot \boldsymbol{n}]=v, \text { on } \Gamma,
\end{equation*}
\begin{equation*}
u=g, \text { on } \partial \Omega.
\end{equation*}
where $\Omega$ is a bounded domain in $R^2$, with Lipschitz boundary $\partial \Omega$, and the interface $\Gamma$ is closed and divides $\Omega$ into two disjoint sub-domains $\Omega^-$ and $\Omega^+$; $w$ and $v$ are two functions
defined only along the interface $\Gamma$. The function $f(\boldsymbol{x},u)$ contains $u$ and denotes the nonlinearity, and has different nonlinear
forms with respect to $u$. $\beta$ is weakly degenerate coefficient functions (degenerate points belong to the interface), it is also mentioned other poor properties such as $\infty \geq \beta \geq 0$ ($\beta$ tends to $0$ on the interface). $[u]=u^{+}(x)-u^{-}(x)=w$ and $[\beta(\boldsymbol{x}) \nabla u \cdot \boldsymbol{n}]=\beta^+(\boldsymbol{x}) \nabla u^+ \cdot \boldsymbol{n}-\beta^-(\boldsymbol{x}) \nabla u^- \cdot \boldsymbol{n}=v$ are the difference of the limiting values of $u(x)$ from $\Omega^+$ and $\Omega^-$ respectively. Finally, $g$ is a determined function on the boundary $\partial \Omega$.


\subsection{DNN-FD method}
In this research, we focus on using DNN to develop fully-decoupled numerical methods for solving degenerate interface problems. First, we divide the domain $\Omega$ into uniform Cartesian meshes, we use DNN to solve examined problems on singular sub-domains near the interface, then extract two decoupled boundary value sub-problems on regular sub-domains with no interface. Those two nonlinear systems can be computed in parallel by finite difference method,
 \begin{figure}
 	\centering
	\includegraphics[width=15cm]{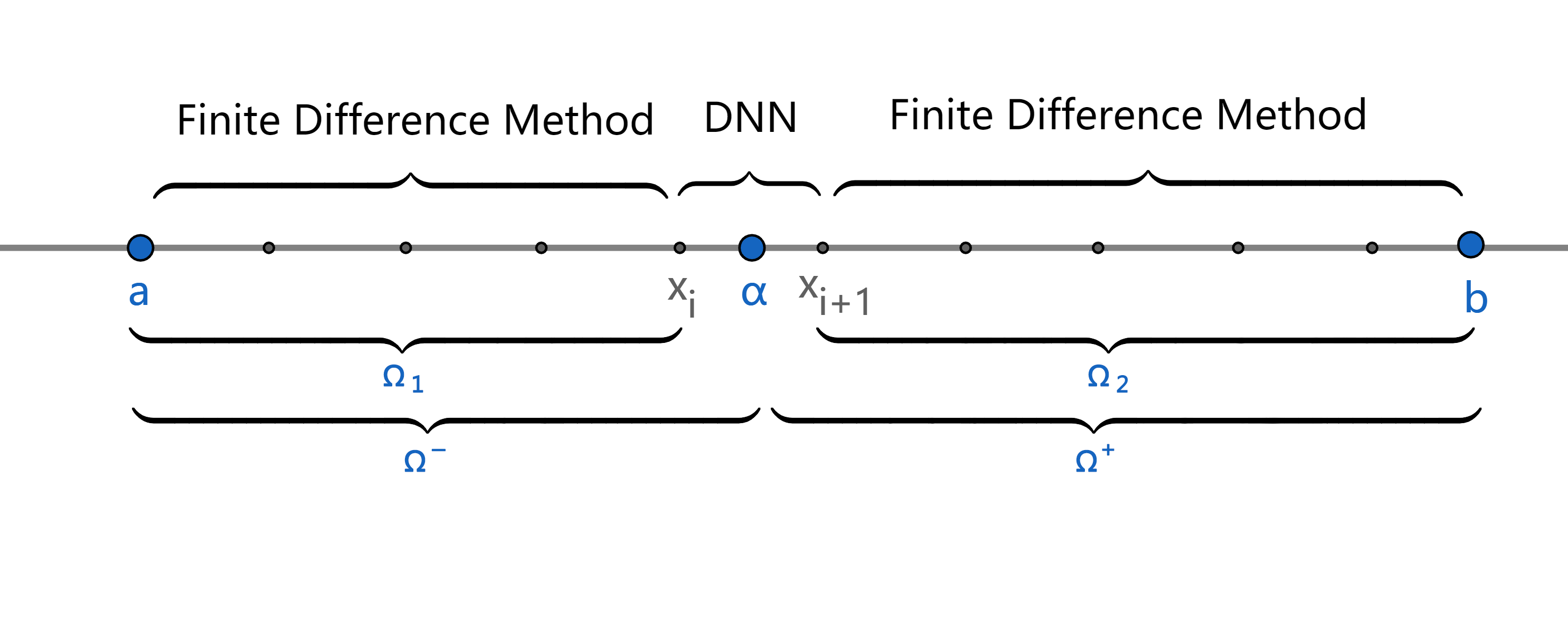}\\
	\caption{A diagram of the method in one-dimensional case}
	\label{1D}
\end{figure}
\begin{equation}
\text { (I) }\left\{\begin{array}{l}
-\nabla \cdot(\beta^-(\boldsymbol{x}) \nabla u^-)=f^-(\boldsymbol{x},u^-), \quad \boldsymbol{x} \in \Omega_1, \\
u^-=u^-_t(\boldsymbol{x};\boldsymbol{\theta}^-), \quad \boldsymbol{x} \in \Gamma^-.
\end{array}\right.
\end{equation}
\begin{equation}
\text { (II) }\left\{\begin{array}{l}
-\nabla \cdot(\beta^+(\boldsymbol{x}) \nabla u^+)=f^+(\boldsymbol{x},u^+), \quad \boldsymbol{x} \in \Omega_2, \\
u^+=u^+_t(\boldsymbol{x};\boldsymbol{\theta}^+), \quad \boldsymbol{x} \in \Gamma^+,\\
u^+=g, \quad \boldsymbol{x} \in \partial \Omega.
\end{array}\right.
\end{equation}
where $f^{\pm}$, $\beta^{\pm}$ and $u^{\pm}$ are the functions in $\Omega^{\pm}$
respectively; $\Omega_{1}$ and $\Omega_{2}$ are regular domains shown in Fig.\ref{1D} and  $u^{\pm}_t(\boldsymbol{x};\boldsymbol{\theta}^{\pm})$ are the result of the deep neural network in the next section.

The proposed method has the advantage of totally decoupling the original problems while using uniform grids. Because our fully decoupled technique is independent of the interface and jump conditions, it not only yields two nondegenerate sub-problems, but it can also easily handle the interface problems with large jump ratios. This method can easily handle both 1D and 2D cases. It is very simple to deal with difficulties like sharp-edge interface issues. While an effective approach for handling general interface problems is not suitable for solving such nonlinear degenerate interface problems, our method can be used robustly and efficiently to both general and degenerate interface problems.


\subsubsection{Deep Neural Network Structure}
In recent years, deep neural network has shown its strong ability in various fields\cite{yu2018deep,handa2016gvnn,pao1989adaptive,collobert2008unified}, mainly reflected in nonlinear fitting ability, high-dimensional data processing ability, excellent fault tolerance ability and strong feature extraction ability. Here, we apply it to the element mesh near the interface to solve the nonlinearity, degeneration and interface singularity of the original problem.

\begin{figure}[htp]
	\centering
	\includegraphics[width=10cm]{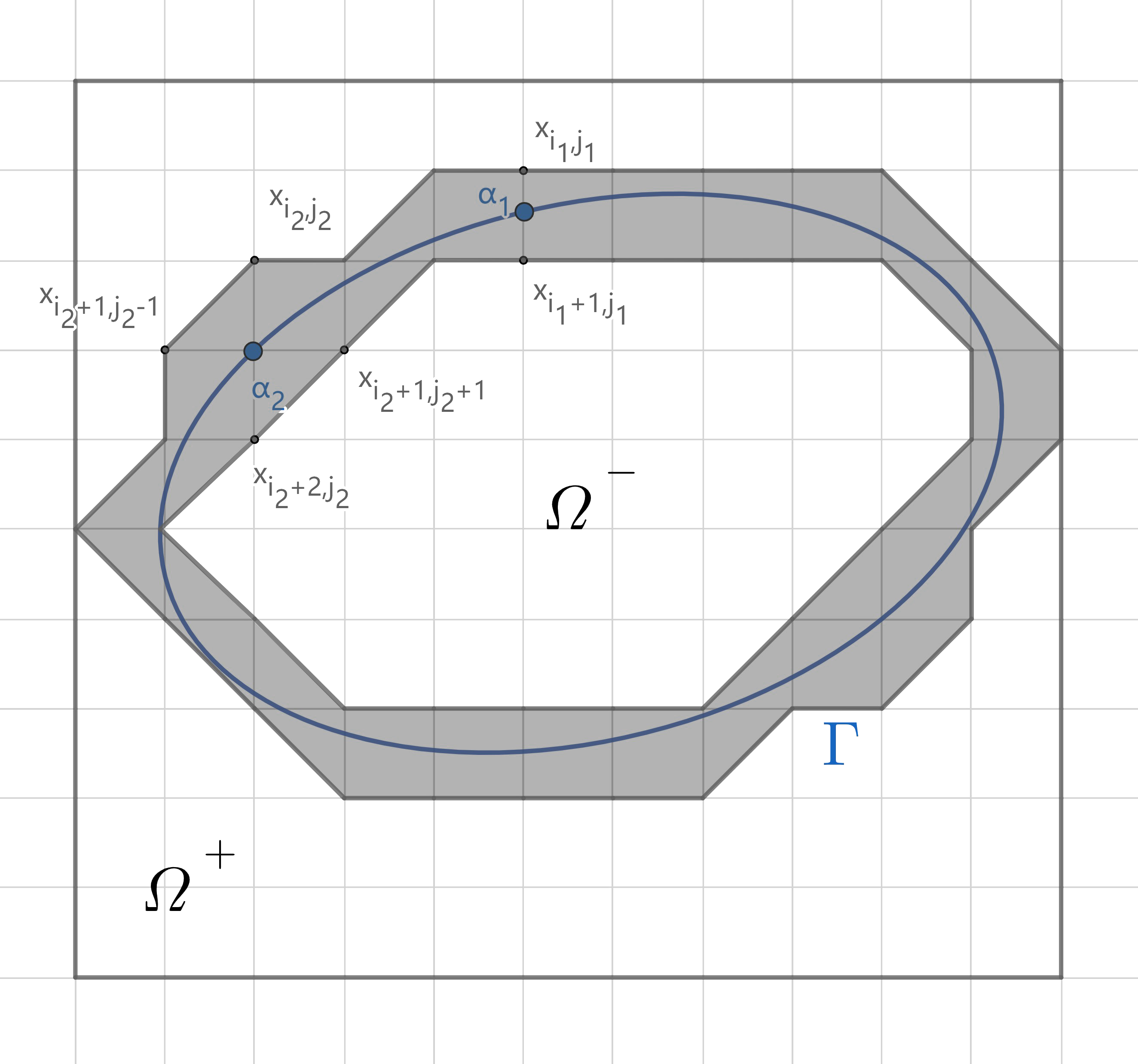}\\
	\caption{A diagram of the method in two-dimensional case}
	\label{2D}
\end{figure}

We apply DNN in the banded degenerate domain composed of near interface element grid in Fig.\ref{2D}. We construct the DNN structure on this domain instead of the whole area to approximate the solution $u$. The reason is that we want to solve the singularity on the interface through the characteristics of DNN, in order to avoid the influence of regular domains on the accuracy of DNN. And the regular domains can be improved by better numerical methods. The problem is naturally separated into two nonsingular sub-problems\cite{lagaris1998artificial,han2017deep,wang2020mesh},
\begin{equation}
u(\boldsymbol{x}) \approx u_t(\boldsymbol{x};\boldsymbol{\theta})= \begin{cases}u^-_t(\boldsymbol{x};\boldsymbol{\theta^{-}}), & \text { if } \boldsymbol{x} \in \Omega^-\setminus \Omega_{1}, \\
u^+_t(\boldsymbol{x};\boldsymbol{\theta^{+}}), & \text { if } \boldsymbol{x} \in \Omega^+\setminus \Omega_{2}.\end{cases}
\end{equation}
\begin{equation}
u_t^{+}\left(\boldsymbol{x} ; \boldsymbol{\theta}^{+}\right)=(|\boldsymbol{x}-\boldsymbol{x}_0|+1)\hat{g}(\boldsymbol{x}_0)+|\boldsymbol{x}-\boldsymbol{x}_0|\hat{u}_t^{+}\left(\boldsymbol{x} ; \boldsymbol{\theta}^{+}\right).
\end{equation}
where $\boldsymbol{\theta}=(\boldsymbol{\theta^-};\boldsymbol{\theta^+}) \in \Theta$, the exact interface is the zero level set of the following level set function $\phi(\boldsymbol{x}_0)=0$. $\hat{g}$ is an extension of $g$ near the interface and $|.|$ is the Euclidean distance. $\hat{u}_t^{+}$ will be obtained from deep learning networks. The construction of equation (3.5) aims to ensure the uniqueness of the solution. Similarly, depending on the shape of the interface, $u^-_t(\boldsymbol{x};\boldsymbol{\theta^{-}})$ will also be constructed correspondingly. If the first jump condition across the interface is homogeneous, only one function $u_t(\boldsymbol{x};\boldsymbol{\theta})$ can be used to approximate the solution $u$.
The structure of DNN with four hidden layers has been given in the Fig.\ref{DNN}. The following is the selection of sampling points, which is divided into two types: one is to select interior points
$\left\{\boldsymbol{x}_{k}\right\}_{k=1}^{M_1}$, $\left\{\boldsymbol{x}_{k}\right\}_{k=1}^{M_2}$
which are random on the degenerate domains; and the other is the nodes $\left\{\boldsymbol{x}_{k}\right\}_{k=1}^{M_3}$
on the element grids. In order to define the discrete loss function, all sampling points
$\left\{\boldsymbol{x}_{k}\right\}_{k=1}^{M_1}$, $\left\{\boldsymbol{x}_{k}\right\}_{k=1}^{M_2}$,
$\left\{\boldsymbol{x}_{k}\right\}_{k=1}^{M_3}$
need to meet the first condition in (\ref{e2}),
\begin{equation}
L_1(\boldsymbol{\theta}):=\frac{1}{M_1+M_3/2} \sum_{k=1}^{M_1+M_3/2}|-
\nabla \cdot \beta^- \nabla u^-_t(\boldsymbol{x}_k;\boldsymbol{\theta})-
f^-(\boldsymbol{x}_k)|^2,\boldsymbol{x} \in \Omega^-\setminus \Omega_{1},
\end{equation}
\begin{equation}
L_2(\boldsymbol{\theta}):=\frac{1}{M_2+M_3/2} \sum_{k=1}^{M_2+M_3/2}|-
\nabla \cdot \beta^+ \nabla u^+_t(\boldsymbol{x}_k;\boldsymbol{\theta})-
f^+(\boldsymbol{x}_k)|^2,\boldsymbol{x} \in \Omega^+\setminus \Omega_{2}.
\end{equation}
The nodes $\left\{\boldsymbol{x}_{k}\right\}_{k=1}^{M_3}$ also need to meet the jump conditions across the interface,
\begin{equation}
L_{3}(\boldsymbol{\theta}):=\frac{2}{M_3} \sum_{k=1}^{M_3}|
u^+_t(\boldsymbol{x}_{i^+_k,j^+_k}; \boldsymbol{\theta})-
u^-_t(\boldsymbol{x}_{i^-_k,j^-_k}; \boldsymbol{\theta})-w|^2,
\end{equation}
\begin{equation}
L_{4}(\boldsymbol{\theta}):=\frac{2}{M_3} \sum_{k=1}^{M_3}|
\beta^+ \nabla u^+_t(\boldsymbol{x}_{i^+_k,j^+_k}; \boldsymbol{\theta})\cdot \boldsymbol{n}-
\beta^- \nabla u^-_t(\boldsymbol{x}_{i^-_n,j^-_n}; \boldsymbol{\theta})\cdot \boldsymbol{n}-v|^2.
\end{equation}
This structure is to solve the singularity and geometric irregularity on the interface. If we sample points directly from the interface, the separated sub-problems will be also degenerate.

In particular, there may be two cases for nodes, the first case is that the intersection of the interface and the grid is not a grid node shown in the Fig.\ref{2D}, such as $\alpha_1$, we can process by nodes close to the intersection in the horizontal or vertical direction,
\begin{equation}
|
u^+_t(\boldsymbol{x}_{i_1,j_1}; \boldsymbol{\theta})-
u^-_t(\boldsymbol{x}_{i_1+1,j_1}; \boldsymbol{\theta})-w|^2,
\end{equation}
\begin{equation}
|
\beta^+ \nabla u^+_t(\boldsymbol{x}_{i_1,j_1}; \boldsymbol{\theta})\cdot \boldsymbol{n}-
\beta^- \nabla u^-_t(\boldsymbol{x}_{i_1+1,j_1}; \boldsymbol{\theta})\cdot \boldsymbol{n}-v|^2.
\end{equation}
The second case is that the interface just intersects with the grid at the node, such as $\alpha_2$. We need to deal with it through the four nodes around it,
\begin{equation}
|
u^+_t(\boldsymbol{x}_{i_2,j_2}; \boldsymbol{\theta})-
u^-_t(\boldsymbol{x}_{i_2+2,j_2}; \boldsymbol{\theta})-w|^2+
u^+_t(\boldsymbol{x}_{i_2+1,j_2-1}; \boldsymbol{\theta})-
u^-_t(\boldsymbol{x}_{i_2+1,j_2+1}; \boldsymbol{\theta})-w|^2,
\end{equation}
\begin{equation}
\begin{aligned}
|
\beta^+ \nabla u^+_t(\boldsymbol{x}_{i_2,j_2}; \boldsymbol{\theta})\cdot \boldsymbol{n}-
\beta^- \nabla u^-_t(\boldsymbol{x}_{i_2+2,j_2}; \boldsymbol{\theta})\cdot \boldsymbol{n}-v|^2+\\
|
\beta^+ \nabla u^+_t(\boldsymbol{x}_{i_2+1,j_2-1}; \boldsymbol{\theta})\cdot \boldsymbol{n}-
\beta^- \nabla u^-_t(\boldsymbol{x}_{i_2+1,j_2+1}; \boldsymbol{\theta})\cdot \boldsymbol{n}-v|^2.
\end{aligned}
\end{equation}
now, we are ready to define the total discrete loss function as follows:
\begin{equation}
L(\boldsymbol{\theta}):=
w_1L_1(\boldsymbol{\theta})+
w_2L_2(\boldsymbol{\theta})+
w_3L_3(\boldsymbol{\theta})+
w_4L_4(\boldsymbol{\theta}),
\end{equation}
where $w_i,i=1,2,3,4$ are weights, which are used to solve the problem with large jump ratios. Therefore, each discrete loss function can be compared by the same order of magnitude. After we get the approximation of the gradient with respect to $\boldsymbol{\theta}_{k}$, we can update each component of $\boldsymbol{\theta}$ as
\begin{equation}
\boldsymbol{\theta}_{k}^{n+1}=\boldsymbol{\theta}_{k}^{n}-\left.\eta \frac{\partial L( \boldsymbol{\theta})}{\partial \boldsymbol{\theta}}\right|_{\boldsymbol{\theta}=\boldsymbol{\theta}^{n}_k},
\end{equation}
where $\boldsymbol{\theta}_{k}$ is any component of $\boldsymbol{\theta}$ and $\eta$ is the learning rate. For the sake of simplicity, $\eta$ is usually taken as $10^{-4}$ unless specified.

\subsubsection{Finite Difference Scheme}
On the regular domain, we can use better numerical methods to improve the accuracy of the whole regions. Here we use the finite difference method\cite{beale2007accuracy}. Take one of these areas as an example,
\begin{equation*}
\text { (II) }\left\{\begin{array}{l}
-\nabla \cdot(\beta^+(\boldsymbol{x}) \nabla u^+)=f^+(\boldsymbol{x},u^+), \quad \boldsymbol{x} \in \Omega_2, \\
u^+=u^+_t(\boldsymbol{x};\boldsymbol{\theta}^+), \quad \boldsymbol{x} \in \Gamma^+,\\
u^+=g, \quad \boldsymbol{x} \in \partial \Omega.
\end{array}\right.
\end{equation*}
Suppose that the function $u^+$ has the following nodes $({x_1}_i,{x_2}_j)$ on the domain $\Omega=[a, b] \times[c, d]$, where
\begin{equation*}
a={x_1}_{0}<{x_1}_{1}<{x_1}_{2}<\cdots<{x_1}_{i}<\cdots<{x_1}_{N-1}<{x_1}_{N}=b,
\end{equation*}
\begin{equation*}
c={x_2}_{0}<{x_2}_{1}<{x_2}_{2}<\cdots<{x_2}_{j}<\cdots<{x_2}_{M-1}<{x_2}_{M}=d.
\end{equation*}
The steps are $h_1$ and $h_2$ respectively, and ${x_1}_{i}={x_1}_{0}+i h_1~(i=0,1, \cdots, N)$, ${x_2}_{j}={x_2}_{0}+j h_2~(j=0,1, \cdots, M)$.
By Taylor formula, numerical calculation usually uses the following first-order central difference quotient and second-order central difference quotient to approximate the first-order partial derivative and second-order partial derivative of the function $u^+$ at the node $({x_1}_i,{x_2}_j)$ respectively,
\begin{equation}
\label{e3}
\delta_{x_1} u^+_{i j}=\frac{u^+_{i+1/2, j}-u^+_{i-1/2, j}}{h_1},
~\delta_{x_2} u^+_{i j}=\frac{u^+_{i, j+1/2}-u^+_{i, j-1/2}}{h_2}.
\end{equation}
\begin{equation}
\label{e4}
\delta_{x_1}^{2} u^+_{i j}=\frac{u^+_{i+1, j}-2 u^+_{i j}+u^+_{i-1, j}}{h_1^{2}},~
\delta_{x_2}^{2} u^+_{i j}=\frac{u^+_{i, j+1}-2 u^+_{i j}+u^+_{i, j-1}}{h_2^{2}}.
\end{equation}
where ${x_1}_{i\pm1/2}={x_1}_{i}\pm {h_1}/2$, ${x_2}_{j\pm1/2}={x_2}_{j}\pm {h_2}/2$, $u^+_{i j}$ is the approximate value of the function $u^+$ at the node.

For the equation (II), the difference quotient is used to approximate the partial derivative at the nodes, and the following difference equations can be obtained on the domain $\Omega_2$:
\begin{equation}
\label{e5}
\delta_{x_1}\left(\beta^+_{i j} \delta_{x_1} u^+_{i j}\right)+\delta_{ x_2}\left(\beta^+_{i j} \delta_{x_2} u^+_{i j}\right)=f^+_{i j},
\end{equation}
where $f^+_{i j}=f^+({x_1}_i,{x_2}_j,u^+_{i j})$. By substituting (\ref{e3}) and (\ref{e4}) into (\ref{e5}), we can get
\begin{equation}
\label{e6}
\begin{aligned}
&\frac{1}{h_1^2}\left(\beta^+_{i+1/2, j} u^+_{i+1, j}-\left(\beta^+_{i+1/2, j}+\beta^+_{i-1/2, j}\right) u^+_{i j}+\beta^+_{i-1/2, j} u^+_{i-1, j}\right)+ \\
&\frac{1}{h_2^{2}}\left(\beta^+_{i, j+1/2} u^+_{i, j+1}-\left(\beta^+_{i, j+1/2}+\beta^+_{i, j-1/2}\right) u^+_{i j}+\beta^+_{i, j-1/2} u^+_{i, j-1}\right)=f^+_{i j}.
\end{aligned}
\end{equation}
where $\beta^+_{i j}=\beta^+ \left({x_1}_{i}, {x_2}_{j}\right)$, $\beta^+_{i \pm 1/2, j}=\beta^+ \left({x_1}_{i \pm 1/2}, {x_2}_{j}\right)$, $\beta^+_{i,j \pm 1/2}=\beta^+ \left({x_1}_{i}, {x_2}_{j \pm 1/2}\right)$, $i=1, \cdots, N-1$, $j=1, \cdots, M-1$.

After discretizing the boundary value conditions, we can get
\begin{equation*}
u^+_{i j}=u^+_t({x_1}_{i}, {x_2}_{j};\boldsymbol{\theta}^+), ~({x_1}_{i}, {x_2}_{j}) \in \left\{\boldsymbol{x}_{k}\right\}_{k=1}^{M_3}.
\end{equation*}
\begin{equation}
u^+_{0 j}=g_{0 j}, u^+_{N j}=g_{N j}, u^+_{i 0}=g_{i 0}, u^+_{i M}=g_{i M},  i=0, \cdots, N, j=0, \cdots, M.
\end{equation}
where $g_{i j}=g({x_1}_i,{x_2}_j)$. Finally, the following iterative method is used to solve (\ref{e6}), set an initial value ${u^+_{i j}}^{(0)}(i=1, \cdots, N-1, j=1, \cdots, M-1)$ and construct the sequence ${u^+_{i j}}^{(m)}(i=1, \cdots, N-1, j=1, \cdots, M-1, m=0,1, \cdots)$ according to the following formula:
\begin{equation}
\begin{aligned}
&\frac{1}{h_1^2}\left(\beta^+_{i+1/2, j} u^{+(m)}_{i+1, j}-\left(\beta^+_{i+1/2, j}+\beta^+_{i-1/2, j}\right) {u^+_{i j}}^{(m)}+\beta^+_{i-1/2, j} u^{+(m)}_{i-1, j}\right)+ \\
&\frac{1}{h_2^{2}}\left(\beta^+_{i, j+1/2} u^{+(m)}_{i, j+1}-\left(\beta^+_{i, j+1/2}+\beta^+_{i, j-1/2}\right) {u^+_{i j}}^{(m)}+\beta^+_{i, j-1/2} u^{+(m)}_{i, j-1}\right)={f^+_{i j}}^{(m)}.
\end{aligned}
\end{equation}


\section{Numerical examples}
\label{s4}
In this section, we present some numerical results to illustrate the expected convergence rates for different configurations. The convergence order of the approximate solutions, as measured by the errors, is denoted by
\begin{equation*}
\text { order }=\log _{2}\left(\left\|u_{2 h}-u\right\|_{L^2} /\left\|u_{h}-u\right\|_{L^2}\right),
\end{equation*}
where $u_{h}$ is the numerical solution with space step size $h$ and $u$ is the analytical solution.

\subsection{1D degenerate interface with homogeneous jump conditions}
\label{example1}
\textbf{Example 4.1.} The degenerate differential equation with the homogeneous interface condition will be solved in $\Omega^{-}=(0,1), \Omega^{+}=(1,2)$, and the interface point $\alpha=1$. The boundary condition and the source function are chosen so that the exact solution is\cite{zhao2021semi}
\begin{equation*}
u(x)=\left\{
\begin{array}{l}
\frac{1}{\tau^{-}}\left(-\exp (1-x)^{1 / 2}+1\right),x \in \Omega^{-}, \\
\frac{1}{\tau^{+}}\left(\exp (x-1)^{1 / 2}-1\right), x \in \Omega^{+}.
\end{array}\right.
\end{equation*}
The coefficient $\beta$ is
\begin{equation*}
\beta=\left\{\begin{array}{l}
\tau^{-}(1-x)^{1 / 2},  x \in \Omega^{-}, \\
\tau^{+}(x-1)^{1 / 2},  x \in \Omega^{+}.
\end{array}\right.
\end{equation*}
Hence, the interface jump conditions,
\begin{equation*}
[u]=w=0, \quad\left[\beta u_{x}\right]=v=0.
\end{equation*}
\begin{table}{}
	\centering
	\caption{$L^2$ errors and convergence orders $\left(\tau^{-} / \tau^{+}=\right.$ $\left.10^{12} / 1\right)$ for Example \ref{example1}.}
	\label{table1}
	\begin{tabular}{lllllll}
		& $\left\|u_{h}-u\right\|_{L^2(\Omega_1)}$ &   & $\left\|u_{h}-u\right\|_{L^2(\Omega_2)}$&  & $\left\|u_{h}-u\right\|_{L^2(\Omega)}$\\
		N & Error & Order & Error & Order & Error & Order \\
		\hline 10 & $4.77 \mathrm{E}-03$ & $-$ & $6.14 \mathrm{E}-03$ & $-$ & $4.39 \mathrm{E}-03$ & $-$ \\
		20 & $1.15 \mathrm{E}-03$ & $2.0508$ & $1.73 \mathrm{E}-03$ & $1.8239$ & $1.22 \mathrm{E}-03$ & $1.8462$ \\
		40 & $3.02 \mathrm{E}-04$ & $1.9286$ & $4.70 \mathrm{E}-04$ & $1.8830$ & $3.10 \mathrm{E}-04$ & $1.9767$ \\
		80 & $7.77 \mathrm{E}-05$ & $1.9608$ & $1.23 \mathrm{E}-04$ & $1.9278$ & $8.07 \mathrm{E}-05$ & $1.9428$ \\
		160 & $1.98 \mathrm{E}-05$ & $1.9668$ & $3.19 \mathrm{E}-05$ & $1.9520$ & $2.03 \mathrm{E}-05$ & $1.9864$
	\end{tabular}
\end{table}
\begin{figure}
	\centering    
	\subfigure[$\tau^{-} / \tau^{+}=10^{12} / 1$]
	{
		\begin{minipage}{6cm}
			\centering
			\includegraphics[width=6cm,clip]{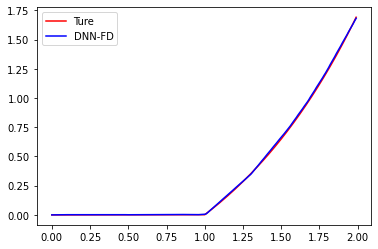}
			
		\end{minipage}
	}
	\subfigure[$\tau^{-} / \tau^{+}=1 / 10^{12}$]
	{
		\begin{minipage}{6cm}
			\centering
			\includegraphics[width=6cm,clip]{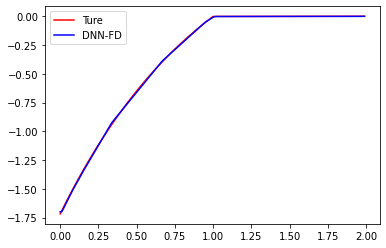}
			
		\end{minipage}
	}
	\vspace{-3mm}
	\caption{Comparison between exact and DNN-FD solutions for Example \ref{example1} (N=160).}
	\label{figure1}
\end{figure}
\begin{table}{}
	\centering
	\caption{$L^2$ errors and convergence orders $\left(\tau^{-} / \tau^{+}=\right.$ $\left.1 / 10^{12}\right)$ for Example \ref{example1}.}
	\label{table2}
	\begin{tabular}{lllllll}
		& $\left\|u_{h}-u\right\|_{L^2(\Omega_1)}$ &   & $\left\|u_{h}-u\right\|_{L^2(\Omega_2)}$&  & $\left\|u_{h}-u\right\|_{L^2(\Omega)}$\\
		N & Error & Order & Error & Order & Error & Order \\
		\hline 10 & $1.42 \mathrm{E}-02$ & $-$ & $7.75 \mathrm{E}-03$ & $-$ & $7.30 \mathrm{E}-03$ & $-$ \\
		20 & $2.94 \mathrm{E}-03$ & $2.2767$ & $2.23 \mathrm{E}-03$ & $1.7977$ & $1.84 \mathrm{E}-03$ & $1.9834$ \\
		40 & $6.61 \mathrm{E}-04$ & $2.1546$ & $5.23 \mathrm{E}-04$ & $2.0903$ & $4.76 \mathrm{E}-04$ & $1.9547$ \\
		80 & $1.57 \mathrm{E}-04$ & $2.0719$ & $1.34 \mathrm{E}-04$ & $1.9623$ & $1.16 \mathrm{E}-04$ & $2.0376$ \\
		160 & $3.60 \mathrm{E}-05$ & $2.1253$ & $3.24 \mathrm{E}-05$ & $2.0512$ & $2.80 \mathrm{E}-05$ & $2.0496$
	\end{tabular}
\end{table}
\begin{figure}
	\centering    
	\subfigure[$\tau^{-} / \tau^{+}=10^{12} / 1$]
	{
		\begin{minipage}{6cm}
			\centering
			\includegraphics[width=6cm,clip]{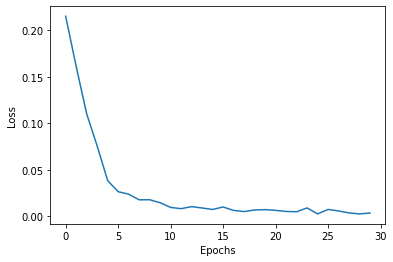}
			
		\end{minipage}
	}
	\subfigure[$\tau^{-} / \tau^{+}=1 / 10^{12}$]
	{
		\begin{minipage}{6cm}
			\centering
			\includegraphics[width=6cm,clip]{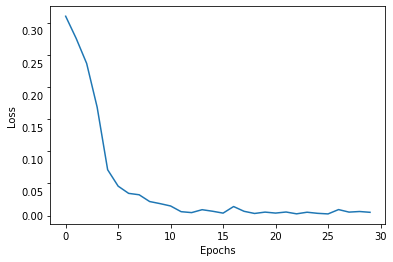}
			
		\end{minipage}
	}
	\vspace{-3mm}
	\caption{The decay of the loss functions for Example \ref{example1}  (N=160).}
	\label{figure2}
\end{figure}
We test the current method for the classical interface problem with homogeneous jump conditions. The network used 4 intermediate layers. The width of each layer is 6 and the number of sampling points is 202, including 200 interior points and two grid nodes. The numerical results of the current method for the very big jump ratios $\left(\tau^{-} / \tau^{+}=\right.$ $\left.10^{12} / 1  ~\text{and}  ~\tau^{-} / \tau^{+}=1 / 10^{12}\right)$ are shown in Table \ref{table1} and Table  \ref{table2} respectively. It can be seen clearly that the convergence orders reach the second order for the numerical solution in $L^2$ norms. Fig.\ref{figure1} shows the comparison between the exact solution and the numerical solution for the very big jump ratios when N=160.
In Fig.\ref{figure2}, we present the decay of the loss function during the training process respectively, eventually the error between the DNN solution and the exact solution reduces to about $O(10^{-4})$ near the interface.

Many other well-known methods usually give the numerical results with the jump ratios $\left(\tau^{-} / \tau^{+}=\right.$ $\left.10^{3} / 1  ~\text{and}  ~\tau^{-} / \tau^{+}=1 / 10^{3}\right)$ for the one-dimensional or two-dimensional interface problems\cite{hou2005numerical}, while it can be calculated by the method used in this paper with the jump ratios $\left(\tau^{-} / \tau^{+}=\right.$ $\left.10^{12} / 1  ~\text{and}  ~\tau^{-} / \tau^{+}=1 / 10^{12}\right)$. The time for the deep neural network required to simulate the function is approximately 1263 seconds when N=160.


\subsection{1D degenerate interface with nonhomogeneous jump conditions}
\label{example2}

\textbf{Example 4.2.} In this example, the computational domain and interface (a point) are the same as in the previous example. The source function $f(x, u)$ are chosen such that the exact solution is as follows\cite{zhao2021semi}:
\begin{equation*}
u(x)=\left\{\begin{array}{l}
u^{-}(x)=\exp \left((1-x)^{2 / 3}\right), x \in \Omega^{-},\\
u^{+}(x)=\exp \left((x-1)^{1 / 2}\right)+5,x \in \Omega^{+}.
\end{array}\right.
\end{equation*}
The coefficient $\beta$ is
\begin{equation*}
\beta=\left\{\begin{array}{l}
\beta^{-}=\tau^{-}(1-x)^{1 / 3}, x \in \Omega^{-} ,\\
\beta^{+}=\tau^{+}(x-1)^{1 / 2}, x \in \Omega^{+}.
\end{array}\right.
\end{equation*}
The experiment satisfies the following jump conditions,
\begin{equation*}
[u]=w=5,~\left[\beta u_{x}\right]=v=\frac{1}{2} \tau^{+}+\frac{2}{3} \tau^{-}.
\end{equation*}
\begin{table}{}
	\centering
	\caption{$L^2$ errors and convergence orders $\left(\tau^{-} / \tau^{+}=\right.$ $\left.10^{12} / 1\right)$ for Example \ref{example2}.}
	\label{table3}
	\begin{tabular}{lllllll}
		& $\left\|u_{h}-u\right\|_{L^2(\Omega_1)}$ &   & $\left\|u_{h}-u\right\|_{L^2(\Omega_2)}$&  & $\left\|u_{h}-u\right\|_{L^2(\Omega)}$\\
		N & Error & Order & Error & Order & Error & Order \\
		\hline 10 & $1.23 \mathrm{E}-02$ & $-$ & $5.80 \mathrm{E}-03$ & $-$ & $6.54 \mathrm{E}-03$ & $-$ \\
		20 & $4.00 \mathrm{E}-03$ & $1.6279$ & $1.57 \mathrm{E}-03$ & $1.8822$ & $1.64 \mathrm{E}-03$ & $1.9962$ \\
		40 & $1.03 \mathrm{E}-03$ & $1.9515$ & $4.07 \mathrm{E}-04$ & $1.9519$ & $4.55 \mathrm{E}-04$ & $1.8471$ \\
		80 & $2.68 \mathrm{E}-04$ & $1.9493$ & $9.93 \mathrm{E}-05$ & $2.0349$ & $1.15 \mathrm{E}-04$ & $1.9853$ \\
		160 & $6.83 \mathrm{E}-05$ & $1.9706$ & $2.61 \mathrm{E}-05$ & $1.9251$ & $3.02 \mathrm{E}-05$ & $1.9303$		
	\end{tabular}
\end{table}
\begin{table}{}
	\centering
	\caption{$L^2$ errors and convergence orders $\left(\tau^{-} / \tau^{+}=\right.$ $\left.1 / 10^{12}\right)$ for Example \ref{example2}.}
	\label{table4}
	\begin{tabular}{lllllll}
	    & $\left\|u_{h}-u\right\|_{L^2(\Omega_1)}$ &   & $\left\|u_{h}-u\right\|_{L^2(\Omega_2)}$&  & $\left\|u_{h}-u\right\|_{L^2(\Omega)}$\\
		N & Error & Order & Error & Order & Error & Order \\
		\hline 10 & $3.15 \mathrm{E}-03$ & $-$ & $1.37 \mathrm{E}-02$ & $-$ & $8.52 \mathrm{E}-03$ & $-$ \\
		20 & $7.70 \mathrm{E}-04$ & $2.0327$ & $4.24 \mathrm{E}-03$ & $1.6935$ & $2.32 \mathrm{E}-03$ & $1.8730$ \\
		40 & $2.22 \mathrm{E}-04$ & $1.7951$ & $9.79 \mathrm{E}-04$ & $2.1154$ & $6.08 \mathrm{E}-04$ & $1.9346$ \\
		80 & $5.74 \mathrm{E}-05$ & $1.9508$ & $2.75 \mathrm{E}-04$ & $1.8321$ & $1.56 \mathrm{E}-04$ & $1.9637$ \\
		160 & $1.45 \mathrm{E}-05$ & $1.9761$ & $6.53 \mathrm{E}-05$ & $2.0745$ & $4.02 \mathrm{E}-05$ & $1.9545$	
	\end{tabular}
\end{table}
\begin{figure}
	\centering    
	\subfigure[$\tau^{-} / \tau^{+}=10^{12} / 1$]
	{
		\begin{minipage}{6cm}
			\centering
			\includegraphics[width=6cm,clip]{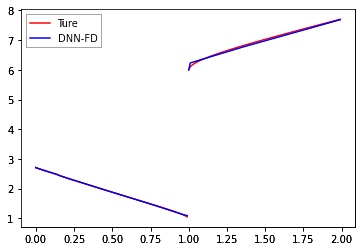}
			
		\end{minipage}
	}
	\subfigure[The decay of the loss functions]
	{
		\begin{minipage}{6cm}
			\centering
			\includegraphics[width=6.5cm,clip]{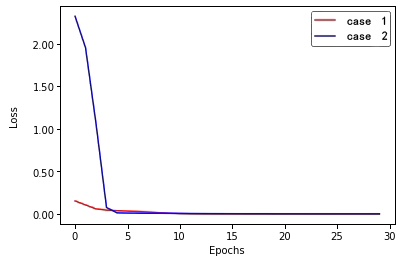}
			
		\end{minipage}
	}
	\vspace{-3mm}
	\caption{Comparison between exact and DNN-FD solutions for Example \ref{example2} (N=80).}
	\label{figure3}
\end{figure}
This is an experiment with nonhomogeneous jump conditions and the requirements for the numerical algorithms problem is higher and stricter to the numerical algorithms. First, we present the convergence order of the variables with large jump ratios $\left(\tau^{-} / \tau^{+}=\right.$ $\left.10^{12} / 1~\text{and}  ~\tau^{-} / \tau^{+}=1 / 10^{12}\right)$ in Table \ref{table3} and Table \ref{table4} namely. It can be seen that the convergence orders for the case of nonhomogeneous jump conditions are the second order. Fig.\ref{figure3}a shows the comparison between the exact solution and the numerical solution for the large jump ratio when N=80. In Fig.\ref{figure3}b, we plot the decay of the $L^2$ norm error between the DNN solution and the exact solution during the training process with the large jump ratio $\left(\tau^{-} / \tau^{+}=\right.$ $\left.10^{12} / 1\right)$ when N=80 (case 2).

Second, to compare with the methods in the literature\cite{zhao2021semi}, we also calculate the results of this experiment with the jump ratio $\left(\tau^{-} / \tau^{+}=\right.$ $\left.10^{7} / 1\right)$. In Fig.\ref{figure3}b, we plot the decay of the loss functions during the training process with jump ratios $\left(\tau^{-} / \tau^{+}=\right.$ $\left.10^{7} / 1~\text{and}  ~\tau^{-} / \tau^{+}=\right.$ $\left.10^{12} / 1\right)$ when N=80. It can be seen that dealing with a smaller jump ratio is more simple and efficient. Finally, using this example, the two methods can calculate homogeneous and nonhomogeneous degenerate problems in one dimension, and the choice of coefficients can be constant, variable, or with singular properties. The advantage of the DNN-FD method is that the jump ratio of the calculated coefficients is bigger than that of the method in\cite{zhao2021semi}. The method can also be extended to two-dimensional degenerate interfaces with the large jump ratio in the next section. This example takes approximately 1298 seconds when N=160, showing that the current method has no essential difference whether the jump conditions are homogeneous or not.


\subsection{2D degenerate interface with nohomogeneous jump conditions}
\label{example3}
\textbf{Example 4.3.} In this example, we consider the interface problem with nonhomogeneous jump conditions. The exact solution is\cite{hou2005numerical}
\begin{equation*}
u(\boldsymbol{x})=\left\{\begin{array}{l}
u^-(\boldsymbol{x})=x_1^{2}+x_2^{2}+2,\boldsymbol{x}\in \Omega^-,\\
u^+(\boldsymbol{x})=1-x_1^{2}-x_2^{2},\boldsymbol{x}\in \Omega^+.
\end{array}\right.
\end{equation*}
The coefficient $\beta$ is
\begin{equation*}
\beta=\left\{\begin{array}{l}
\beta^{-}=\tau^{-}(-\cos(x_{1}^{2}+x_{2}^{2}-(0.5)^{2})+1),\boldsymbol{x}\in \Omega^-;\\
\beta^{+}=\tau^{+}(3-x_1x_2),\boldsymbol{x}\in \Omega^+.
\end{array}\right.
\end{equation*}
where $\Omega^-=\{\boldsymbol{x}|| \boldsymbol{x} \mid<0.5\}, \Omega^+=\Omega \backslash \Omega^{-}, \Omega=[-1,1] \times[-1,1]$, and $r=\sqrt{x_{1}^{2}+x_{2}^{2}}$. The exact interface is the zero level set of the following level set function,
\begin{equation*}
\phi(\boldsymbol{x})=x_{1}^{2}+x_{2}^{2}-(0.5)^{2}.
\end{equation*}
\begin{table}{}
	\centering
	\caption{$L^2$ errors and convergence orders $\left(\tau^{-} / \tau^{+}=\right.$ $\left.10^{10} / 1\right)$ for Example \ref{example3}.}
	\label{table5}
	\begin{tabular}{lllllll}
		& $\left\|u_{h}-u\right\|_{L^2(\Omega_1)}$ &   & $\left\|u_{h}-u\right\|_{L^2(\Omega_2)}$&  & $\left\|u_{h}-u\right\|_{L^2(\Omega)}$\\
		N	& Error & Order & Error & Order & Error & Order \\
		\hline 10 $\times$ 10& $1.45 \mathrm{E}-02$ & $-$ & $1.15 \mathrm{E}-02$ & $-$ & $8.16 \mathrm{E}-03$ & $-$ \\
		20 $\times$ 20 & $4.19 \mathrm{E}-03$ & $1.7978$ & $3.32 \mathrm{E}-03$ & $1.7981$ & $2.45 \mathrm{E}-03$ & $1.7335$ \\
		40 $\times$ 40 & $1.09 \mathrm{E}-03$ & $1.9396$ & $8.56 \mathrm{E}-04$ & $1.9553$ & $6.29 \mathrm{E}-04$ & $1.9636$ \\
		80 $\times$ 80 & $2.83 \mathrm{E}-04$ & $1.9458$ & $2.23 \mathrm{E}-04$ & $1.9366$ & $1.61 \mathrm{E}-04$ & $1.9598$ \\
		160 $\times$ 160& $7.40 \mathrm{E}-05$ & $1.9367$ & $5.73 \mathrm{E}-05$ & $1.9651$ & $4.15 \mathrm{E}-05$ & $1.9598$
	\end{tabular}
\end{table}
\begin{figure}
	\centering    
	\subfigure[DNN-FD Solution]
	{
	\begin{minipage}{5cm}
		\centering
		\includegraphics[width=5cm,clip]{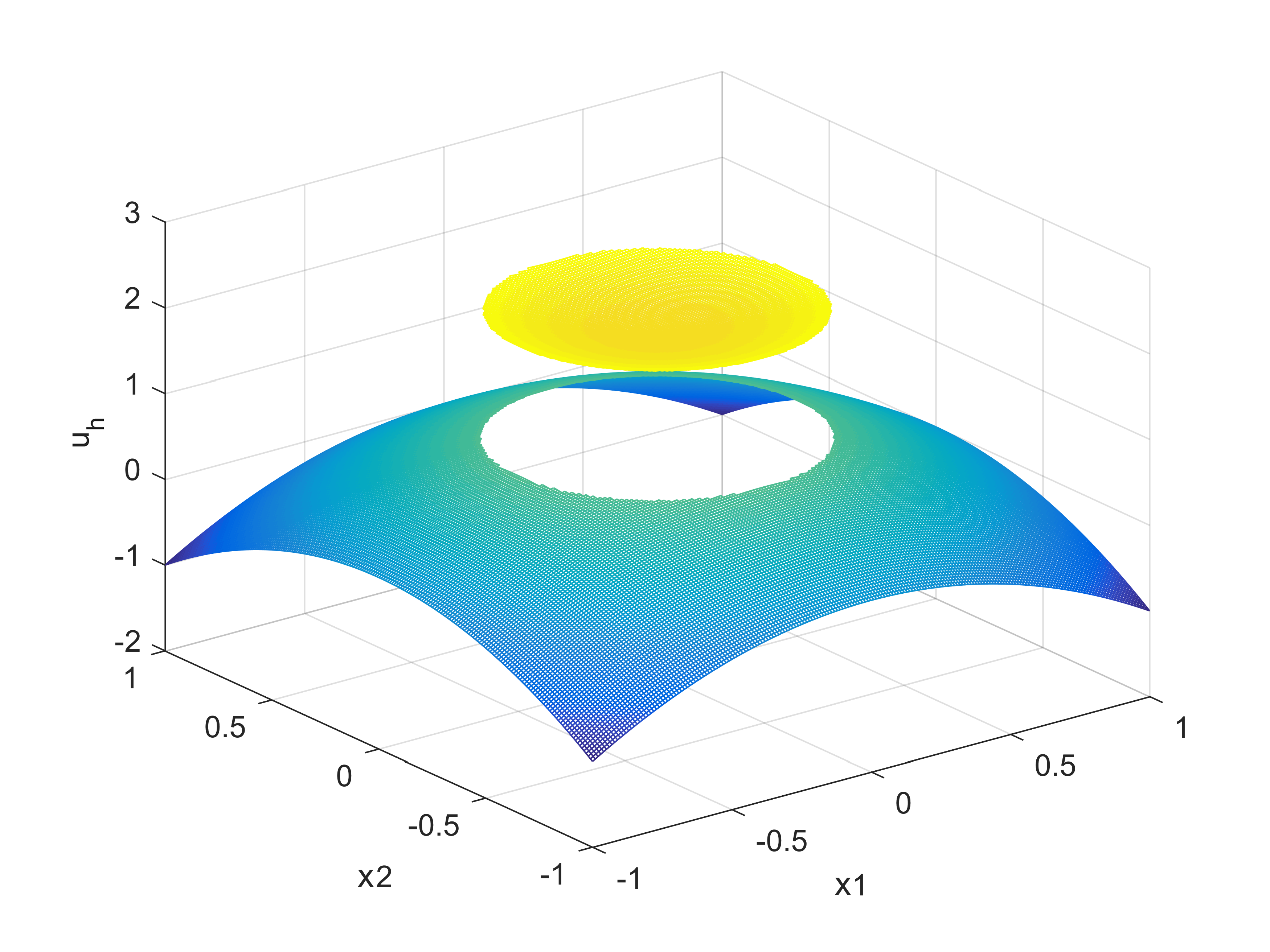}
		
	\end{minipage}
}
    \subfigure[Exact Solution]
   {
	\begin{minipage}{5cm}
		\centering
		\includegraphics[width=5cm,clip]{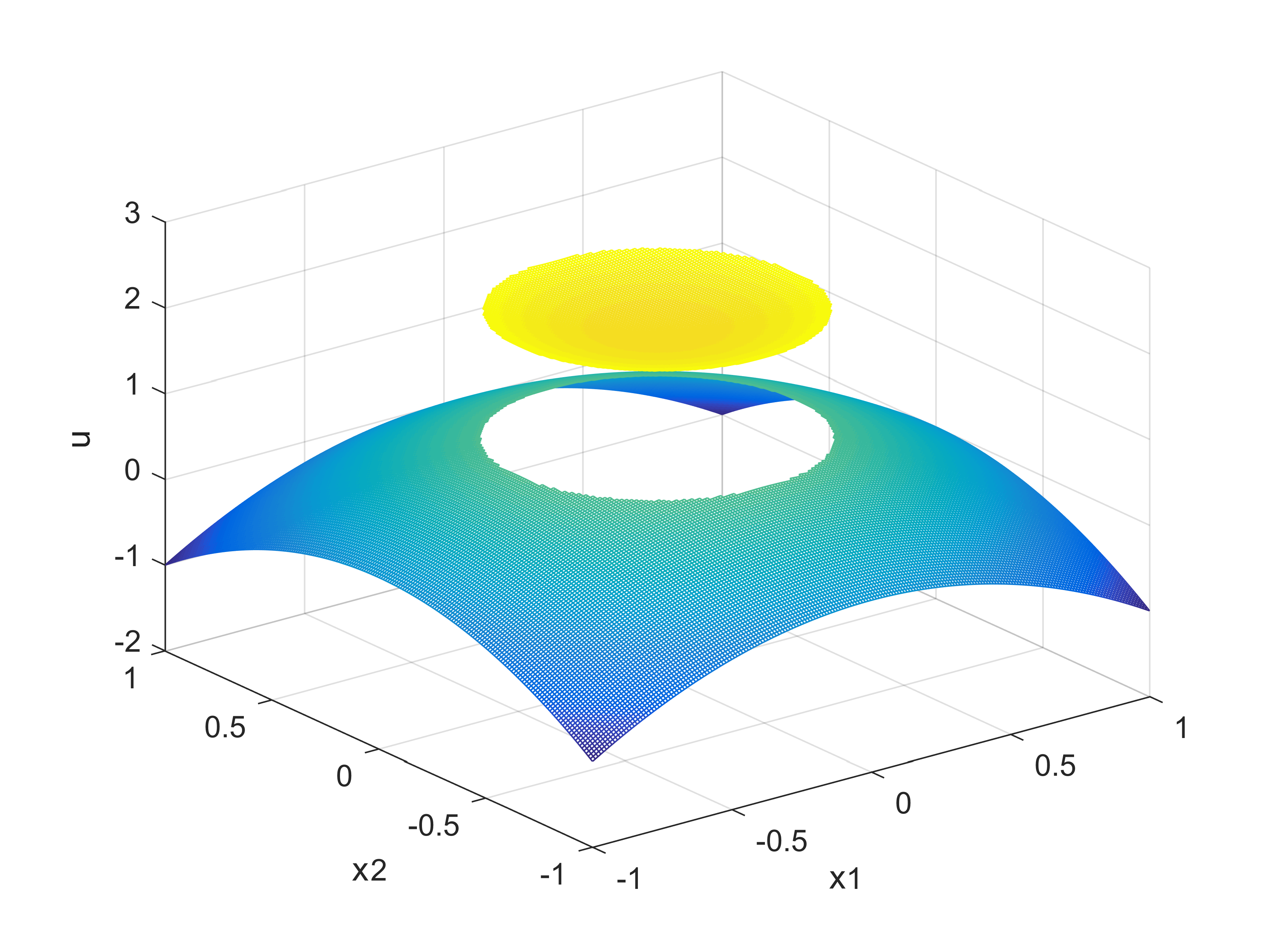}
		
	\end{minipage}
}
	\vspace{-3mm}
	\caption{Comparison between exact and DNN-FD solutions for Example \ref{example3} when N=160 $\left(\tau^{-} / \tau^{+}=1 / 10^{10}\right)$.}
	\label{figure7}
\end{figure}
\begin{table}{}
	\centering
	\caption{$L^2$ errors and convergence orders $\left(\tau^{-} / \tau^{+}=\right.$ $\left.1 / 10^{10}\right)$ for Example \ref{example3}.}
	\label{table6}
	\begin{tabular}{lllllll}
		& $\left\|u_{h}-u\right\|_{L^2(\Omega_1)}$ &   & $\left\|u_{h}-u\right\|_{L^2(\Omega_2)}$&  & $\left\|u_{h}-u\right\|_{L^2(\Omega)}$\\
		N	& Error & Order & Error & Order & Error & Order \\
		\hline 10 $\times$ 10 & $2.45 \mathrm{E}-02$ & $-$ & $1.37 \mathrm{E}-02$ & $-$ & $2.52 \mathrm{E}-02$ & $-$ \\
		20 $\times$ 20 & $6.12 \mathrm{E}-03$ & $2.0020$ & $3.70 \mathrm{E}-03$ & $1.8898$ & $6.22 \mathrm{E}-03$ & $2.0208$ \\
		40 $\times$ 40 & $1.53 \mathrm{E}-03$ & $2.0000$ & $9.55 \mathrm{E}-04$ & $1.9564$ & $1.54 \mathrm{E}-03$ & $2.0114$ \\
		80 $\times$ 80 & $3.82 \mathrm{E}-04$ & $2.0001$ & $2.52 \mathrm{E}-04$ & $1.9202$ & $3.84 \mathrm{E}-04$ & $2.0059$ \\
		160 $\times$ 160 & $9.57 \mathrm{E}-05$ & $2.0000$ & $5.92 \mathrm{E}-05$ & $2.0905$ & $9.59 \mathrm{E}-05$ & $2.0028$
	\end{tabular}
\end{table}
\begin{figure}
	\centering    
	\subfigure[$\tau^{-} / \tau^{+}=10^{10} / 1$]
	{
		\begin{minipage}{6cm}
			\centering
			\includegraphics[width=6cm,clip]{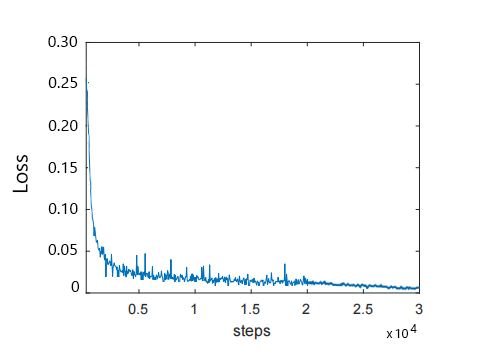}
			
		\end{minipage}
	}
	\subfigure[$\tau^{-} / \tau^{+}=1 / 10^{10}$]
	{
		\begin{minipage}{6cm}
			\centering
			\includegraphics[width=6cm,clip]{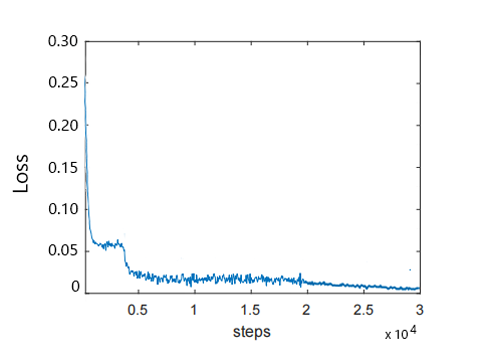}
			
		\end{minipage}
	}
	\vspace{-3mm}
	\caption{The decay of the loss functions for Example \ref{example3} (N=160).}
	\label{figure4}
\end{figure}
We reconstruct the example from the literature\cite{hou2005numerical} to degenerate it near the interface. It is a two-dimensional degenerate elliptic equation with nonhomogeneous jump conditions. The network used 6 intermediate layers. The width of each layer is 15 and the number of sampling interior points is 2000. In the running of the SGD method, we generate a new batch every 10 steps of updating. The numerical results of the present method for large jump ratios $\left(\tau^{-} / \tau^{+}=\right.$ $\left.10^{10} / 1 ~\text{and}~ \tau^{-} / \tau^{+}=1 / 10^{10}\right)$ are shown in Table \ref{table5} and Table \ref{table6} respectively. It can be seen that the convergence orders for the case of nonhomogeneous jump conditions are the second order. Fig.\ref{figure7} shows the comparison between the exact solution and the numerical solution for the large jump ratio $\left(\tau^{-} / \tau^{+}=1 / 10^{10}\right)$ when N=160. In Fig.\ref{figure4}, we plot the decay of the loss functions during the training process with large jump ratios $\left(\tau^{-} / \tau^{+}=\right.$ $\left.10^{10} / 1~\text{and}~\tau^{-} / \tau^{+}=1 / 10^{10}\right)$ when N=160. The two-dimensional case is more difficult than the one-dimensional case and takes more sampling points, but there is no essential difference in methods. The error between the DNN solution and the exact solution is also reduced to approximately $O(10^{-4})$ near the interface. This example shows that this method can be effectively extended to two-dimensional or even higher dimensional degenerate interface problems, and can also effectively solve the coefficients with the large jump ratio.


\subsection{2D nondegenerate interface with homogeneous jump conditions}
\label{example4}
\textbf{Example 4.4.} In this example, we consider the nondegenerate interface problem with high contrast diffusion coefficients with homogeneous jump conditions. The exact solution is\cite{he2022mesh}
\begin{equation*}
u(\boldsymbol{x})=\left\{\begin{array}{l}
u^-(\boldsymbol{x})=\frac{r^{3}}{\beta^-}, \boldsymbol{x} \in \Omega^-, \\
u^+(\boldsymbol{x})=\frac{r^{3}}{\beta^+}+\left(\frac{1}{\beta^-}-\frac{1}{\beta^+}\right) (0.5)^{3}, \boldsymbol{x} \in \Omega^+.
\end{array}\right.
\end{equation*}
where $\Omega^-=\{\boldsymbol{x}|| \boldsymbol{x} \mid<0.5\}, \Omega^+=\Omega \backslash \Omega^-, \Omega=[-1,1] \times[-1,1]$, and $r=\sqrt{x_{1}^{2}+x_{2}^{2}}$. The exact interface is the zero level set of the following level set function,
\begin{equation*}
\phi(\boldsymbol{x})=x_{1}^{2}+x_{2}^{2}-(0.5)^{2}.
\end{equation*}
\begin{table}{}
	\centering
	\caption{$L^2$ errors and convergence orders $\left(\beta^{-} / \beta^{+}=\right.$ $\left.10^{10} / 1\right)$ for Example \ref{example4}.}
	\label{table7}
	\begin{tabular}{lllllll}
		& $\left\|u_{h}-u\right\|_{L^2(\Omega_1)}$ &   & $\left\|u_{h}-u\right\|_{L^2(\Omega_2)}$&  & $\left\|u_{h}-u\right\|_{L^2(\Omega)}$\\
		N & Error & Order & Error & Order & Error & Order \\
		\hline 10 $\times$ 10 & $1.44 \mathrm{E}-02$ & $-$ & $4.36 \mathrm{E}-02$ & $-$ & $1.37 \mathrm{E}-03$ & $-$ \\
		20 $\times$ 20 & $3.56 \mathrm{E}-03$ & $2.0221$ & $1.15 \mathrm{E}-03$ & $1.9252$ & $3.27 \mathrm{E}-03$ & $2.0646$ \\
		40 $\times$ 40 & $8.98 \mathrm{E}-04$ & $1.9891$ & $2.96 \mathrm{E}-04$ & $1.9548$ & $8.26 \mathrm{E}-04$ & $1.9857$ \\
		80 $\times$ 80 & $2.26 \mathrm{E}-04$ & $1.9891$ & $2.08 \mathrm{E}-05$ & $1.9727$ & $7.55 \mathrm{E}-04$ & $1.9849$ \\
		160 $\times$ 160 & $5.68 \mathrm{E}-05$ & $1.9926$ & $1.91 \mathrm{E}-05$ & $1.9840$ & $5.27 \mathrm{E}-05$ & $1.9869$
	\end{tabular}
\end{table}
\begin{figure}
	\centering    
	\subfigure[DNN-FD Solution]
	{
		\begin{minipage}{5cm}
			\centering
			\includegraphics[width=5cm,clip]{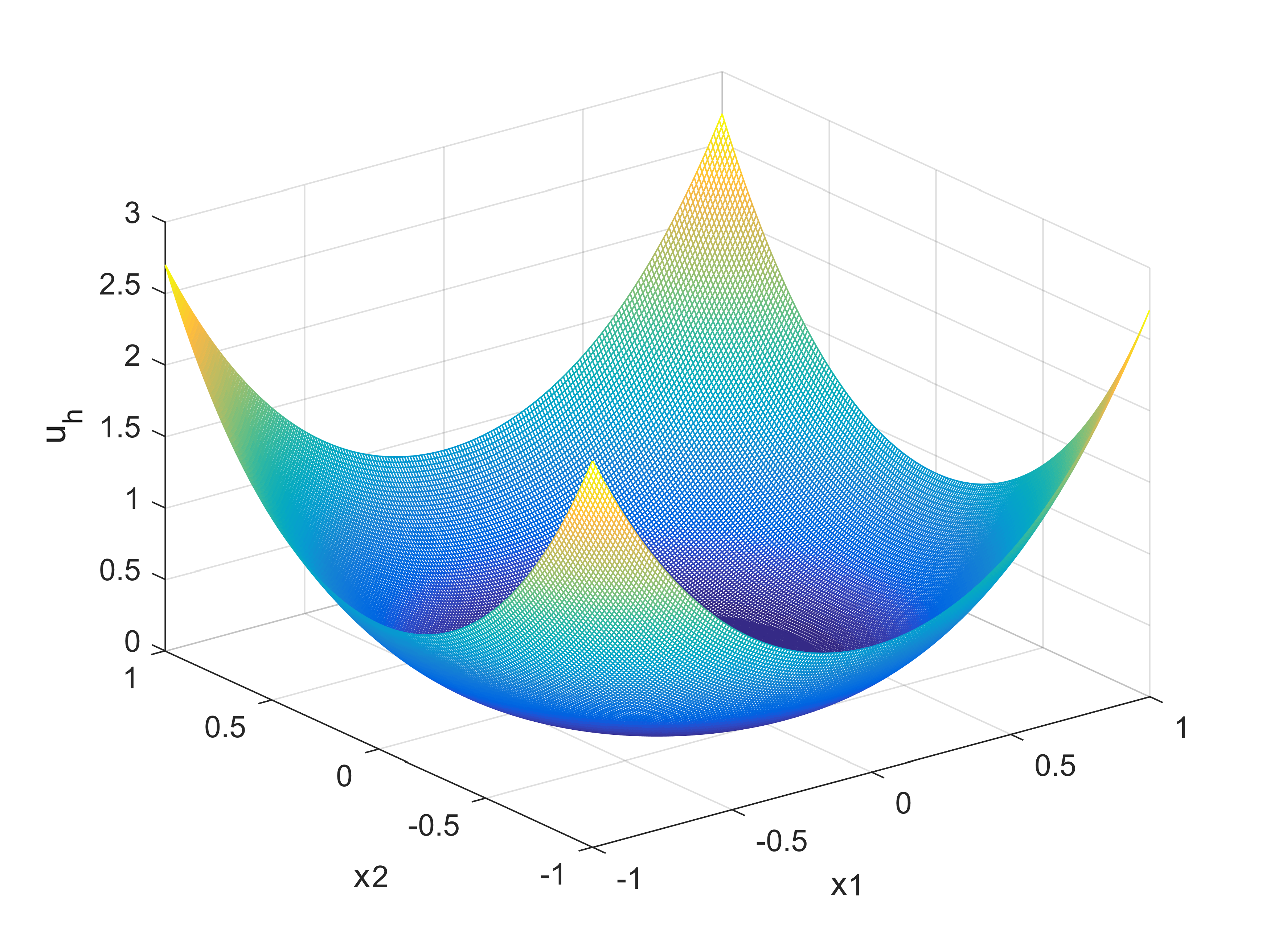}
			
		\end{minipage}
	}
	\subfigure[Exact Solution]
	{
		\begin{minipage}{5cm}
			\centering
			\includegraphics[width=5cm,clip]{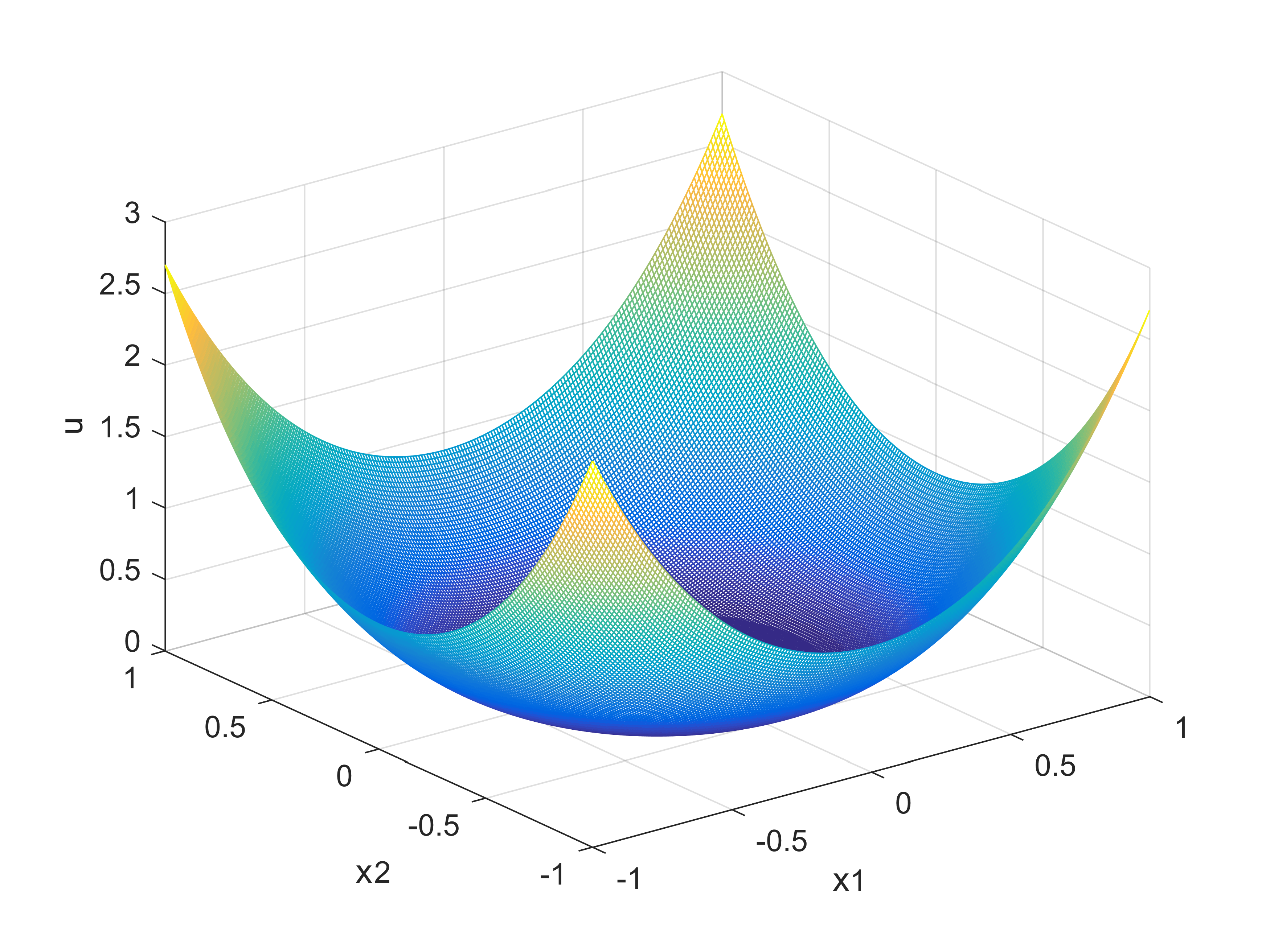}
			
		\end{minipage}
	}
	\vspace{-3mm}
	\caption{Comparison between exact and DNN-FD solutions for Example \ref{example4} when N=160 $\left(\beta^{-} / \beta^{+}=\right.$ $\left.10^{10} / 1\right)$.}
	\label{figure10}
\end{figure}
\begin{table}{}
	\centering
	\caption{$L^2$ errors and convergence orders $ \left(\beta^{-} / \beta^{+}=\right.$ $\left.1 / 10^{10}\right)$ for Example \ref{example4}.}
	\label{table8}
	\begin{tabular}{lllllll}
		& $\left\|u_{h}-u\right\|_{L^2(\Omega_1)}$ &   & $\left\|u_{h}-u\right\|_{L^2(\Omega_2)}$&  & $\left\|u_{h}-u\right\|_{L^2(\Omega)}$\\
		N & Error & Order & Error & Order & Error & Order \\
		\hline 10 $\times$ 10 & $2.12 \mathrm{E}-02$ & $-$ & $1.10 \mathrm{E}-03$ & $-$ & $1.43 \mathrm{E}-02$ & $-$ \\
		20 $\times$ 20 & $5.47 \mathrm{E}-03$ & $1.9531$ & $2.74 \mathrm{E}-03$ & $2.0102$ & $3.59 \mathrm{E}-03$ & $1.9984$ \\
		40 $\times$ 40 & $1.33 \mathrm{E}-03$ & $2.0349$ & $6.91 \mathrm{E}-04$ & $1.9912$ & $8.96 \mathrm{E}-04$ & $2.0032$ \\
		80 $\times$ 80 & $3.26 \mathrm{E}-04$ & $2.0310$ & $1.73 \mathrm{E}-05$ & $1.9975$ & $2.35 \mathrm{E}-04$ & $1.9304$ \\
		160 $\times$ 160 & $7.91 \mathrm{E}-05$ & $2.0471$ & $4.30 \mathrm{E}-05$ & $2.0075$ & $6.14 \mathrm{E}-05$ & $1.9364$
	\end{tabular}
\end{table}
\begin{figure}
	\centering    
	\subfigure[DNN-FD Solution]
	{
		\begin{minipage}{5cm}
			\centering
			\includegraphics[width=5cm,clip]{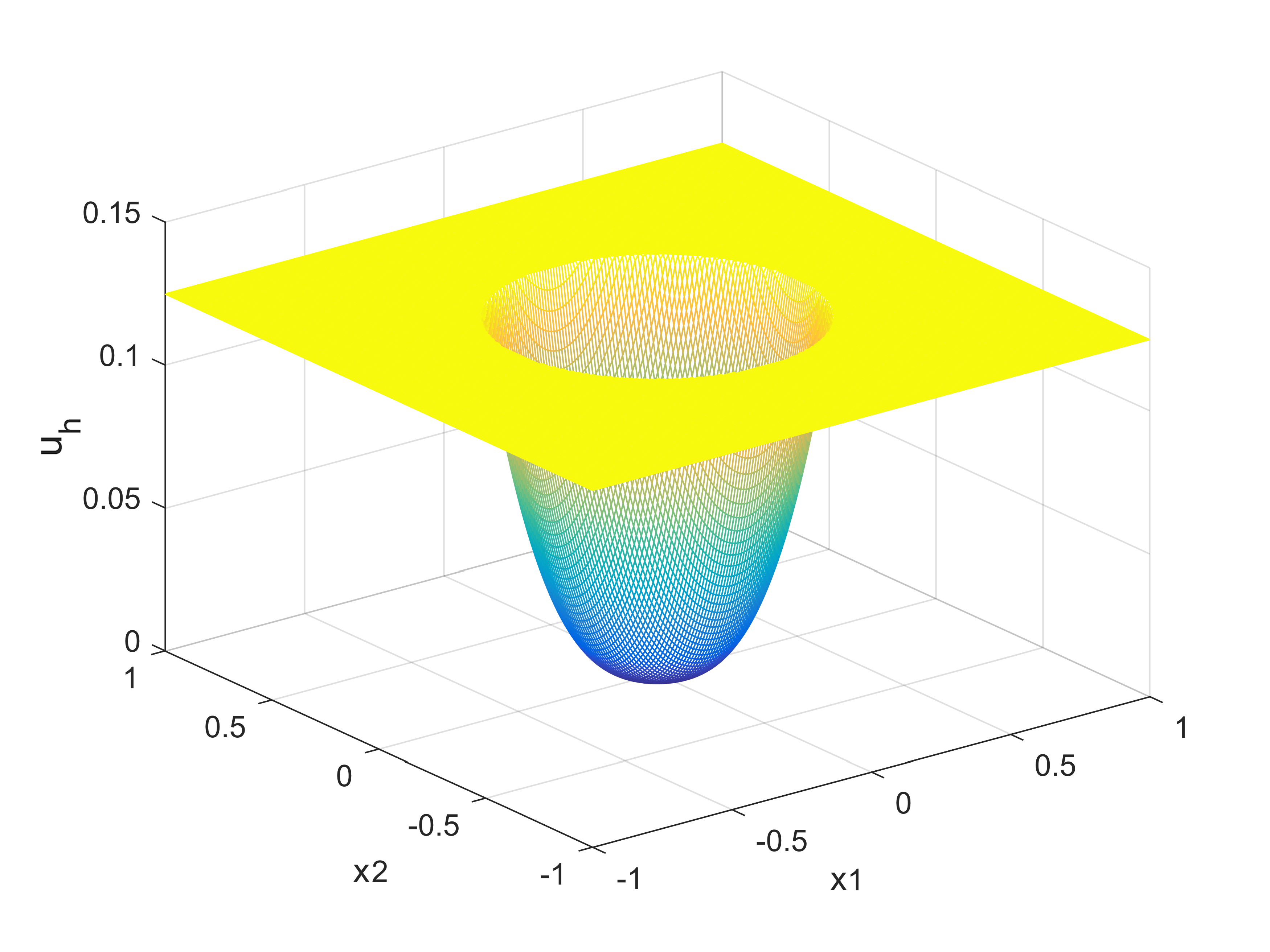}
			
		\end{minipage}
	}
	\subfigure[Exact Solution]
	{
		\begin{minipage}{5cm}
			\centering
			\includegraphics[width=5cm,clip]{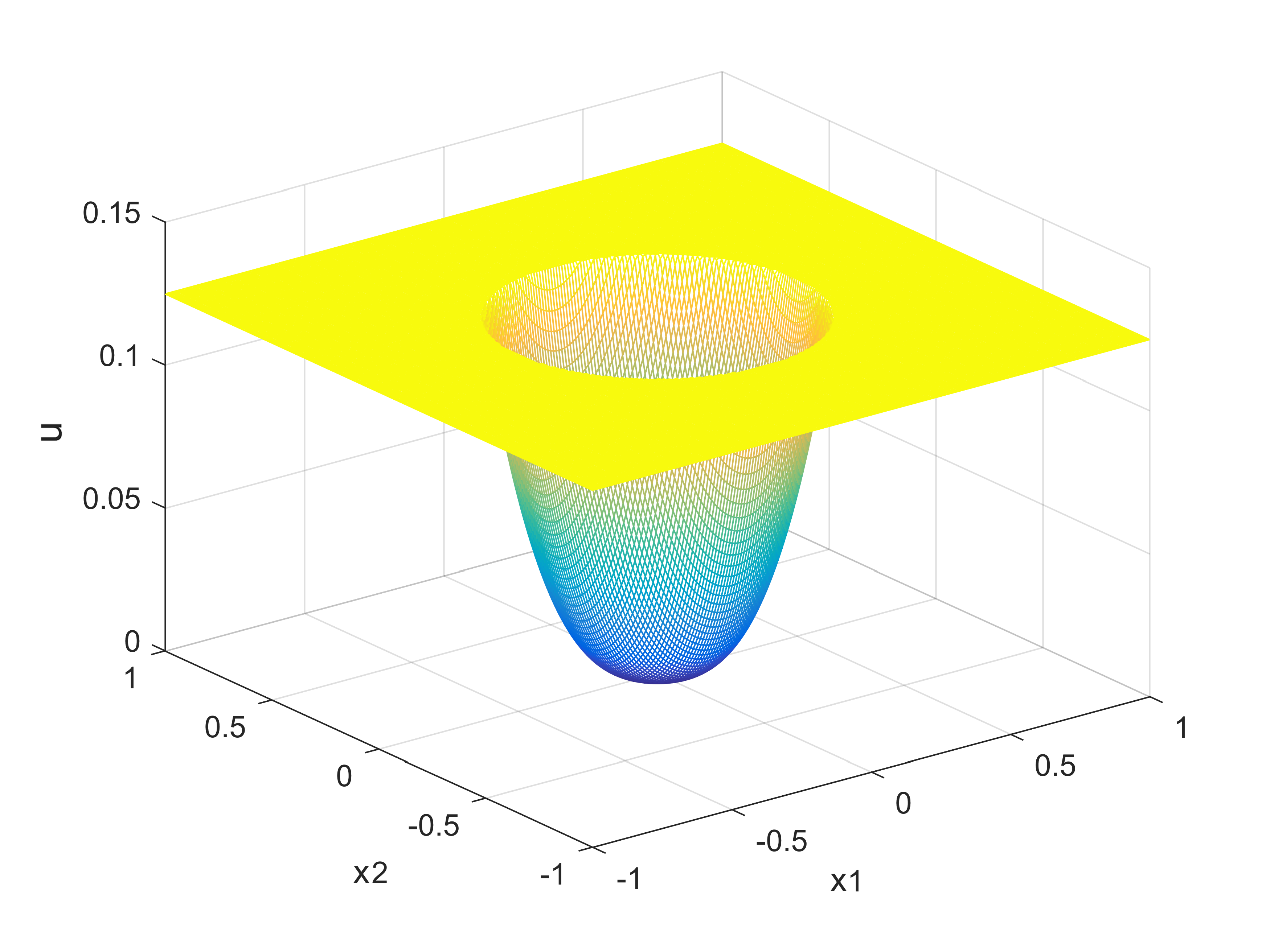}
			
		\end{minipage}
	}
	\vspace{-3mm}
	\caption{Comparisons between exact and DNN-FD solutions for Example \ref{example4} when N=160 $\left(\beta^{-} / \beta^{+}=\right.$ $\left.1 / 10^{10}\right)$.}
	\label{figure11}
\end{figure}
The method used in this paper can compute not only degenerate problems, but also nondegenerate problems. The numerical results of the present method for large jump ratios $\left(\beta^{-} / \beta^{+}=\right.$ $\left.10^{10} / 1~\text{and}~\beta^{-} / \beta^{+}=1 / 10^{10}\right)$ are shown in Table \ref{table7} and Table \ref{table8} respectively. It can be seen easily that the numerical solution has second-order convergence in the $L^2$ norm.

Fig.\ref{figure10} and Fig.\ref{figure11} show the comparison between the exact solution and the numerical solution for large jump ratios $\left(\tau^{-} / \tau^{+}=\right.$ $\left.10^{10} / 1\right)$ and $\left(\tau^{-} / \tau^{+}=1 / 10^{10}\right)$ when N=160 respectively. Due to the application of numerical methods on regular domains, the accuracy of this method is higher than that in \cite{he2022mesh}, and because of the fully decoupled format, it can handle the problem with higher coefficients and the larger jump ratio.


\subsection{2D nondegenerate flower shape interface}
\label{example5}
\textbf{Example 4.5.} In this example, we consider the flower shape interface problem. The exact solution is\cite{hou2005numerical}
\begin{equation*}
u(\boldsymbol{x})=\left\{\begin{array}{l}
u^-(\boldsymbol{x})=7 x_1^{2}+7 x_2^{2}+6,\boldsymbol{x}\in \Omega^-, \\
u^+(\boldsymbol{x})=5-5 x_1^{2}-5 x_2^{2},\boldsymbol{x}\in \Omega^+.
\end{array}\right.
\end{equation*}
The coefficient $\beta$ is
\begin{equation*}
\beta=\left\{\begin{array}{l}
\beta^{-}=\left(x_1^{2}-x_2^{2}+3\right) / 7,\boldsymbol{x}\in \Omega^-,\\
\beta^{+}=(x_1 x_2+2) / 5,\boldsymbol{x}\in \Omega^+.
\end{array}\right.
\end{equation*}
The exact interface is the zero level set of the following level set function,
\begin{equation*}
\begin{aligned}
&\phi=(x_1-0.02 \sqrt{5})^{2}+(x_2-0.02 \sqrt{5})^{2}-(0.5+0.2 \sin (5 \theta))^{2}, \\
&\text { with }\left\{\begin{array}{l}
x(\theta)=0.02 \sqrt{5}+(0.5+0.2 \sin (5 \theta)) \cos (\theta), \\
y(\theta)=0.02 \sqrt{5}+(0.5+0.2 \sin (5 \theta)) \sin (\theta),
\end{array} \quad \theta \in[0,2 \pi).\right.
\end{aligned}
\end{equation*}
\begin{table}{}
	\centering
	\caption{$L^2$ errors and convergence orders of the DNN-FD method for Example \ref{example5}.}
	\label{table9}
	\begin{tabular}{lllllll}
		& $\left\|u_{h}-u\right\|_{L^2(\Omega_1)}$ &   & $\left\|u_{h}-u\right\|_{L^2(\Omega_2)}$&  & $\left\|u_{h}-u\right\|_{L^2(\Omega)}$\\
		N & Error & Order & Error & Order & Error & Order \\
		\hline 10 $\times$ 10 & $1.03 \mathrm{E}-02$ & $-$ & $1.36 \mathrm{E}-02$ & $-$ & $1.30 \mathrm{E}-02$ & $-$ \\
		20 $\times$ 20 & $2.63 \mathrm{E}-03$ & $1.9754$ & $3.61 \mathrm{E}-03$ & $1.9145$ & $3.51 \mathrm{E}-03$ & $1.8912$ \\
		40 $\times$ 40 & $6.56 \mathrm{E}-04$ & $1.9834$ & $8.85 \mathrm{E}-04$ & $2.0285$ & $8.98 \mathrm{E}-04$ & $1.9652$ \\
		80 $\times$ 80 & $1.66 \mathrm{E}-04$ & $2.0000$ & $2.19 \mathrm{E}-04$ & $2.0098$ & $2.32 \mathrm{E}-04$ & $1.9522$ \\
		160 $\times$ 160 & $4.16 \mathrm{E}-05$ & $1.9979$ & $5.18 \mathrm{E}-05$ & $2.0847$ & $5.96 \mathrm{E}-05$ & $1.9609$
	\end{tabular}
\end{table}
\begin{figure}
	\centering    
	\subfigure[DNN-FD Solution]
	{
		\begin{minipage}{5cm}
			\centering
			\includegraphics[width=5cm,clip]{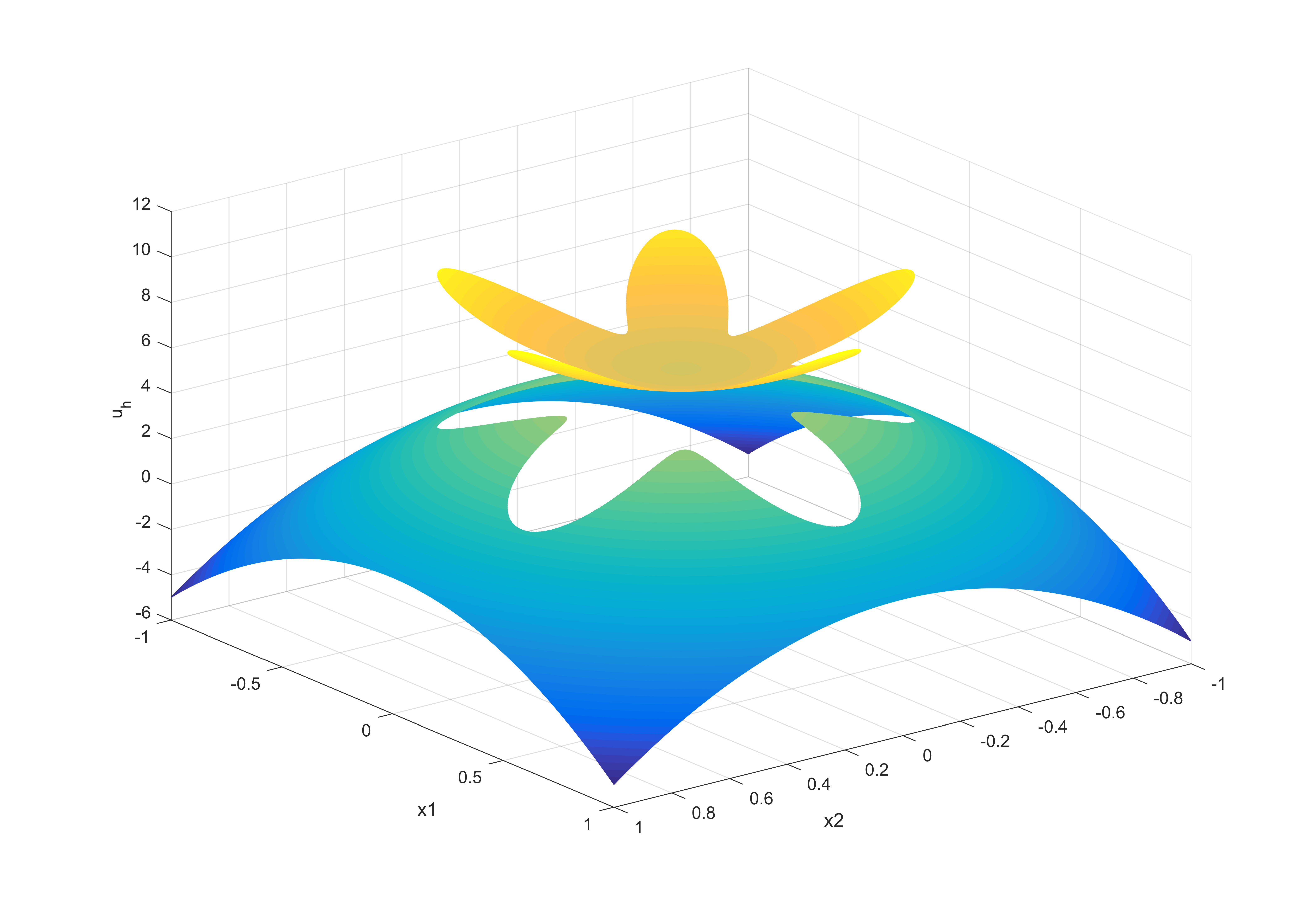}
			
		\end{minipage}
	}
	\subfigure[Exact Solution]
	{
		\begin{minipage}{5cm}
			\centering
			\includegraphics[width=5cm,clip]{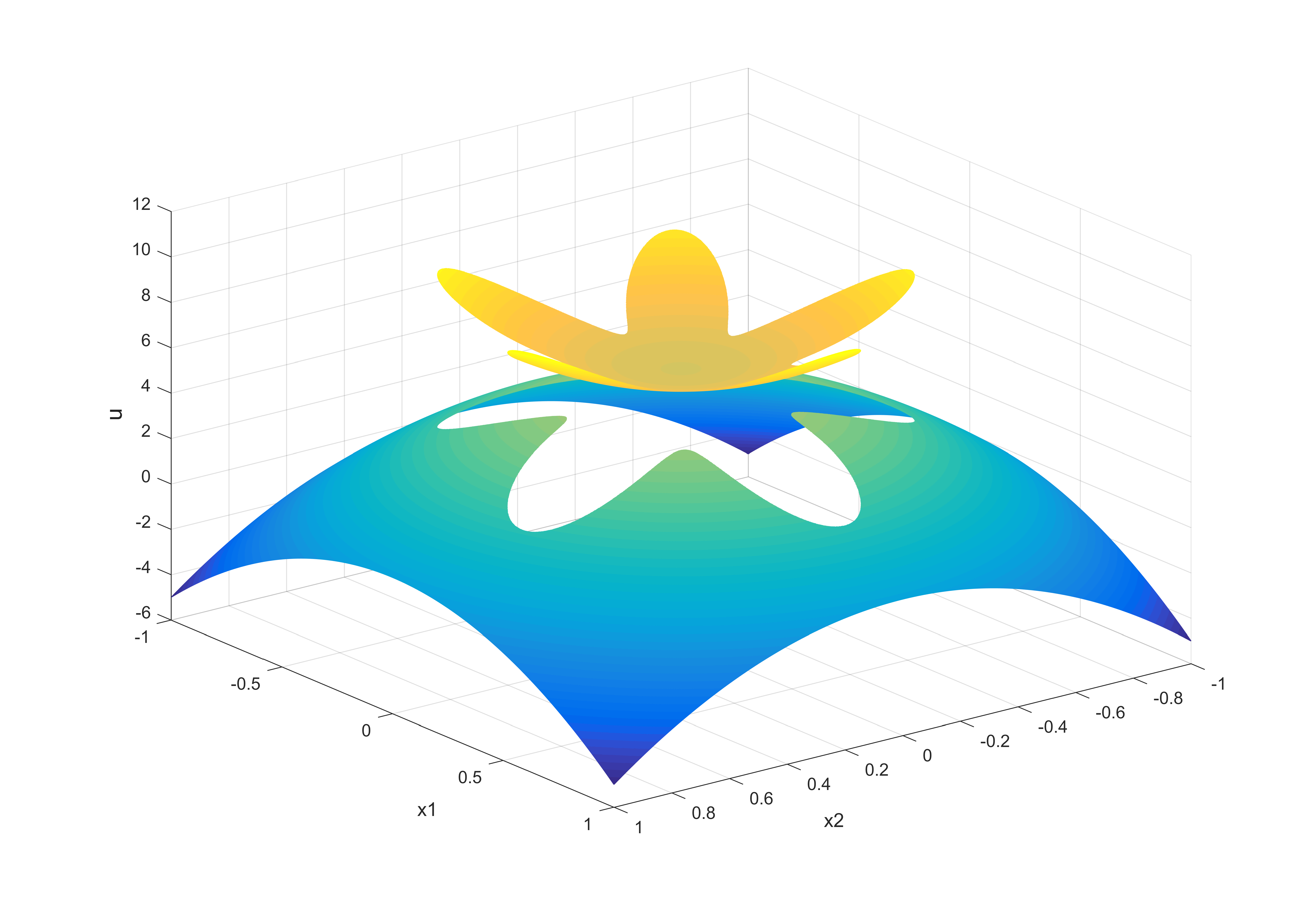}
			
		\end{minipage}
	}
	\vspace{-3mm}
	\caption{Comparison between exact and DNN-FD solutions for Example \ref{example5} when N=160.}
	\label{figure12}
\end{figure}
\begin{figure}
	\centering    
	\subfigure[Sampling points in $\Omega^-$]
	{
		\begin{minipage}{5cm}
			\centering
			\includegraphics[width=5cm,clip]{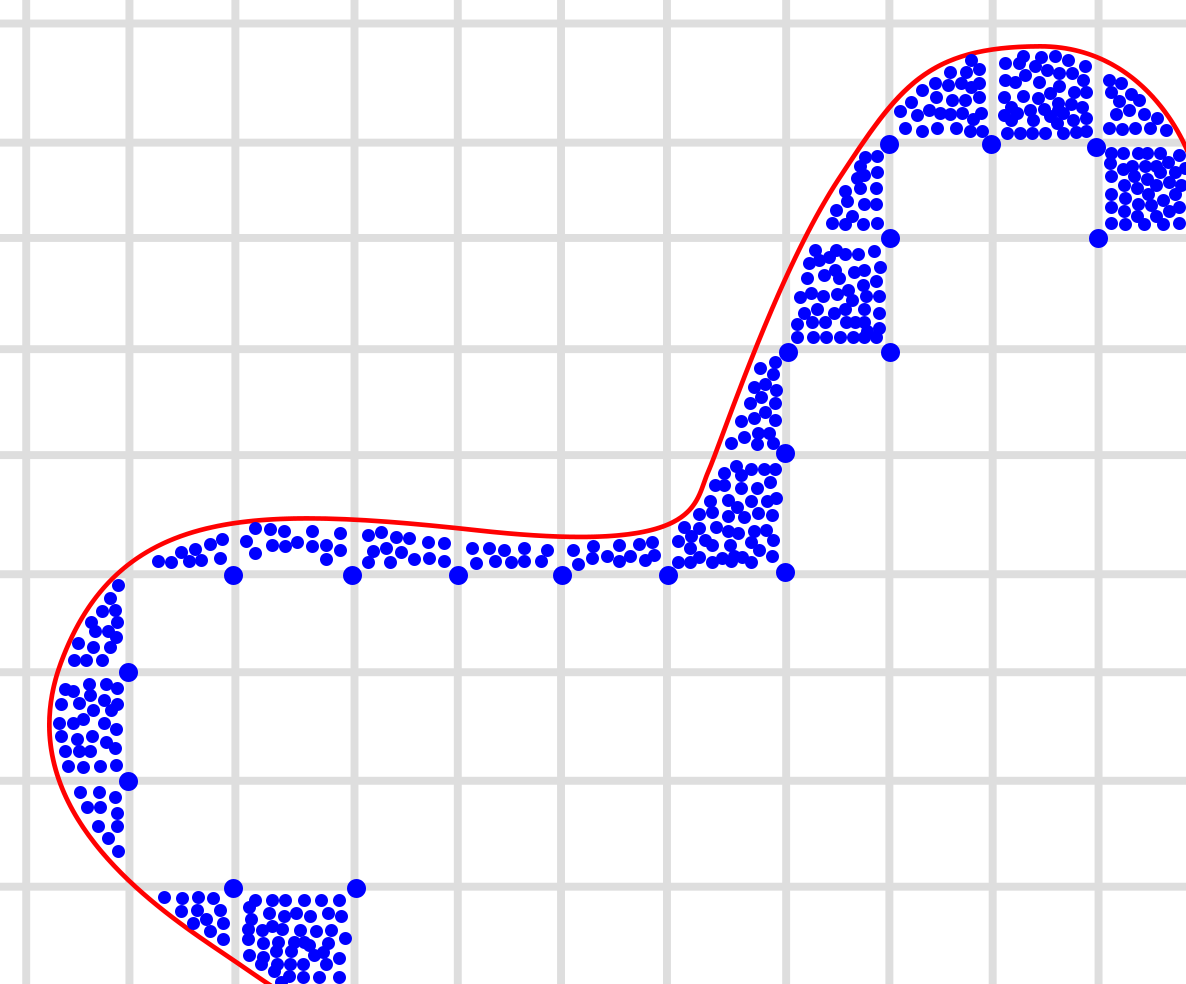}
			
		\end{minipage}
	}
	\subfigure[Sampling points in $\Omega^+$]
	{
		\begin{minipage}{5cm}
			\centering
			\includegraphics[width=5cm,clip]{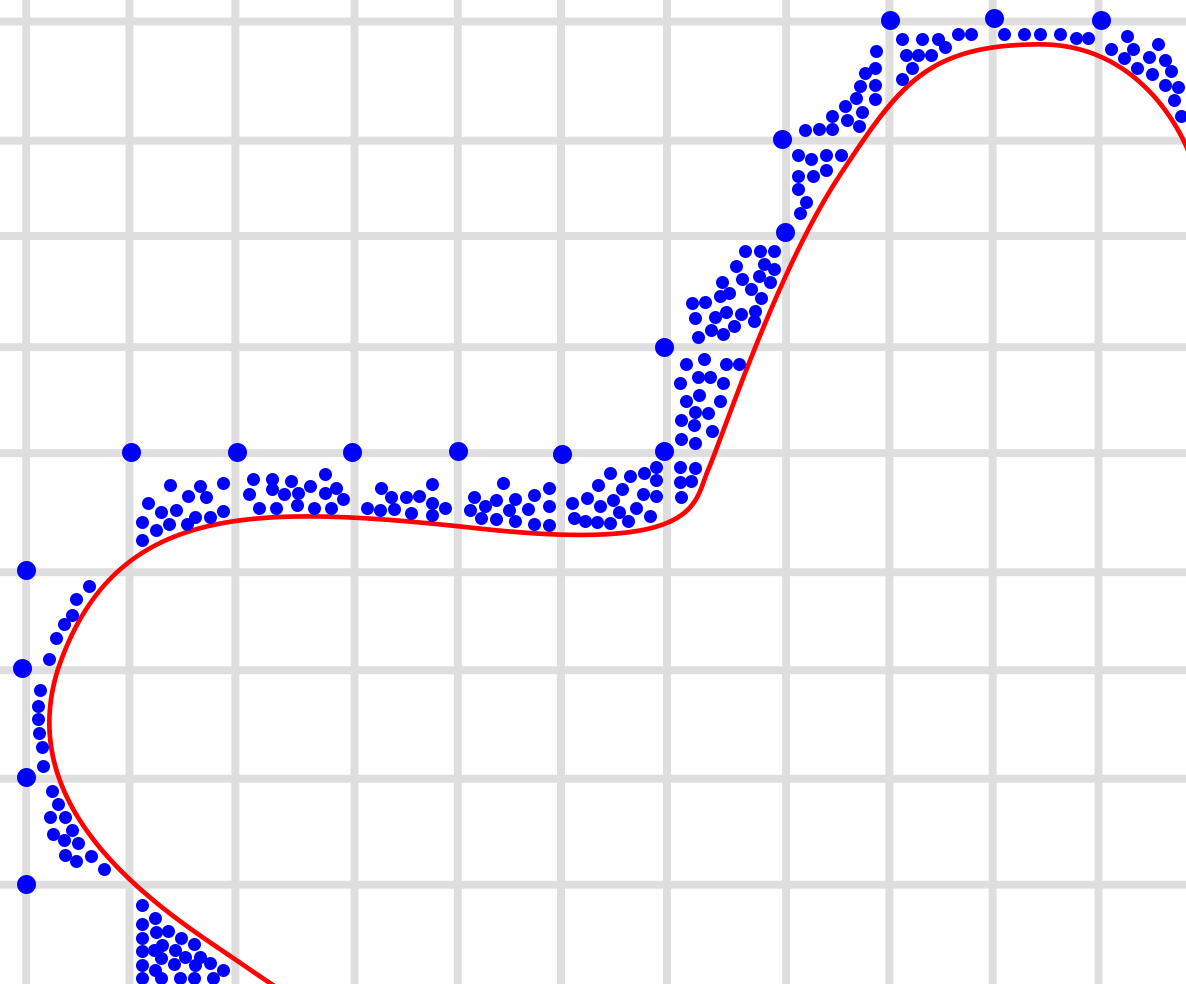}
			
		\end{minipage}
	}
	\vspace{-3mm}
	\caption{The diagram of sampling points for Example \ref{example5}.}
	\label{figure13}
\end{figure}
The peculiarity of this example is that the problem has a complex smooth interface. It is designed to examine the performance of the DNN-FD method in dealing with geometric irregularities. Our method also has advantages in dealing with complex interface problems. This method becomes simple and efficient by applying a deep neural network near the interface. We present a grid refinement analysis in Table \ref{table9} that successfully reached the second order. Fig.\ref{figure13} shows the sampling points in the area of the method in this paper. It can be seen from the figure that we will set more sampling points near the curve with the large radian. Similarly, as dealing with the singularity and non-smoothness of the interface, we will set more sampling points. We take the points by sections based on different degeneracies, large jump ratios, and other conditions to show the properties of the interface well. Fig.\ref{figure12} shows the comparison between the exact solution and the numerical solution when N=160.


\subsection{2D nondegenerate happy-face interface}
\label{example9}
\textbf{Example 4.6.} In this example, we consider the following more general self-adjoint elliptic interface problem,
\begin{equation*}
-\nabla \cdot(\beta(\boldsymbol{x}) \nabla u(\boldsymbol{x}))+\sigma(\boldsymbol{x}) u(\boldsymbol{x})=f(\boldsymbol{x}) \text {, in } \Omega.
\end{equation*}
The example is a happy-face interface and the coefficients $\beta^{\pm}$ are symmetric positive definite matrices. The exact solution is\cite{wang2021bilinear,hou2005numerical}
\begin{equation*}
u(\boldsymbol{x})=\left\{\begin{array}{l}
u^-(\boldsymbol{x})=7 x_1^{2}+7 x_2^{2}+1,\boldsymbol{x}\in \Omega^-,\\
u^+(\boldsymbol{x})=5-5 x_1^{2}-5 x_2^{2},\boldsymbol{x}\in \Omega^+.
\end{array}\right.
\end{equation*}
The coefficient $\beta$ is
\begin{equation*}
\beta^{+}(\boldsymbol{x})=\left(\begin{array}{ll}
x_{1} x_{2}+2 & x_{1} x_{2}+1 \\
x_{1} x_{2}+1 & x_{1} x_{2}+3
\end{array}\right), \beta^{-}(\boldsymbol{x})=\left(\begin{array}{ll}
x_{1}^{2}-x_{2}^{2}+3 & x_{1}^{2}-x_{2}^{2}+1 \\
x_{1}^{2}-x_{2}^{2}+1 & x_{1}^{2}-x_{2}^{2}+4
\end{array}\right).
\end{equation*}
The exact interface can be viewed in the literature\cite{hou2005numerical}.
The other coefficient $\sigma$ is
\begin{equation*}
\sigma(\boldsymbol{x})=\left\{\begin{array}{l}
\sigma^-(\boldsymbol{x})=x_{1} x^{2}+1, \boldsymbol{x}\in \Omega^-, \\
\sigma^+(\boldsymbol{x})=x_{1}^{2}+x_{2}^{2}-2,\boldsymbol{x}\in \Omega^+.
\end{array}\right.
\end{equation*}
\begin{table}{}
	\centering
	\caption{$L^2$ errors and convergence orders of the DNN-FD method for Example \ref{example9}.}
	\label{table14}
	\begin{tabular}{lllllll}
		& $\left\|u_{h}-u\right\|_{L^2(\Omega_1)}$ &   & $\left\|u_{h}-u\right\|_{L^2(\Omega_2)}$&  & $\left\|u_{h}-u\right\|_{L^2(\Omega)}$\\
		N & Error & Order & Error & Order & Error & Order \\
		\hline 10 $\times$ 10 & $1.50 \mathrm{E}-02$ & $-$ & $7.13 \mathrm{E}-02$ & $-$ & $2.37 \mathrm{E}-02$ & $-$ \\
		20 $\times$ 20 & $3.63 \mathrm{E}-03$ & $2.0449$ & $1.31 \mathrm{E}-03$ & $2.4440$ & $7.46 \mathrm{E}-03$ & $1.6659$ \\
		40 $\times$ 40 & $1.42 \mathrm{E}-03$ & $1.3528$ & $4.75 \mathrm{E}-04$ & $1.4647$ & $1.36 \mathrm{E}-03$ & $2.4552$ \\
		80 $\times$ 80 & $3.64 \mathrm{E}-04$ & $1.9644$ & $1.21 \mathrm{E}-04$ & $1.9729$ & $3.23 \mathrm{E}-04$ & $2.0752$ \\
		160 $\times$ 160 & $9.15 \mathrm{E}-05$ & $1.9942$ & $2.97 \mathrm{E}-05$ & $2.0255$ & $8.36 \mathrm{E}-05$ & $1.9503$
	\end{tabular}
\end{table}
\begin{figure}
	\centering    
	\subfigure[DNN-FD Solution]
	{
		\begin{minipage}{5cm}
			\centering
			\includegraphics[width=5cm,clip]{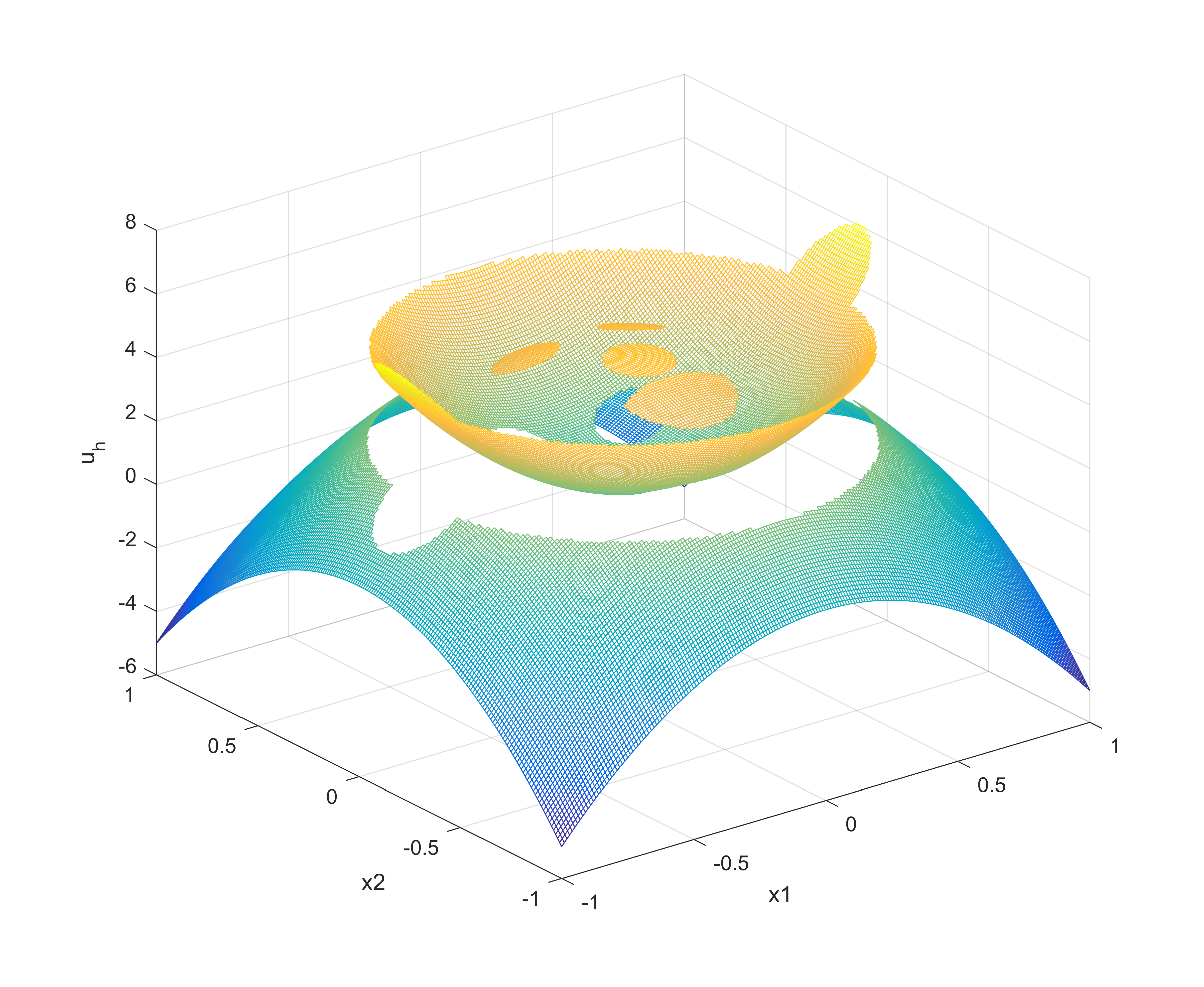}
			
		\end{minipage}
	}
	\subfigure[Exact Solution]
	{
		\begin{minipage}{5cm}
			\centering
			\includegraphics[width=5cm,clip]{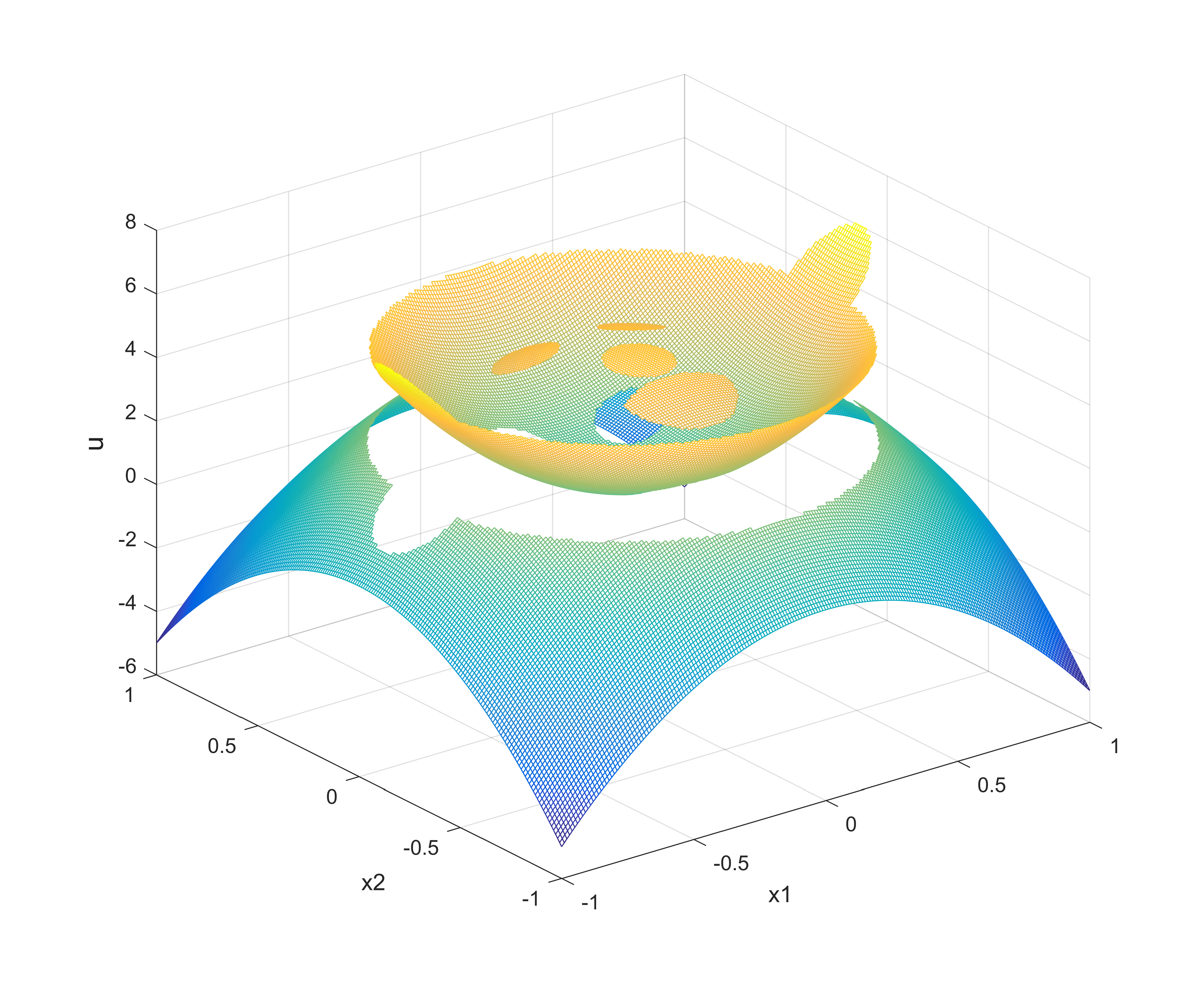}
			
		\end{minipage}
	}
	\vspace{-3mm}
	\caption{Comparison between exact and DNN-FD solutions for Example \ref{example9} when N=160.}
	\label{figure17}
\end{figure}
The difficulty of the example is that the interfaces have kinks around ears and mouth. We present the convergence results in Table \ref{table14}. Numerical results indicate that the DNN-FD solution always converges to the exact solution with second-order accuracy. And the exact solution and the numerical solution are compared in Fig.\ref{figure17} when N=160.


\subsection{2D nondegenerate sharp-edged interface}
\label{example6}
\textbf{Example 4.7.} In this example, we consider the nonsmooth interface problem. The exact solution is\cite{wang2021new,hou2010numerical}
\begin{equation*}
u(\boldsymbol{x})=\left\{\begin{array}{l}
u^-(\boldsymbol{x})=7 x_1^{2}+7 x_2^{2}+6, \\
u^+(\boldsymbol{x})=\begin{cases}x_{1}+x_{2}+1, & \text { if } x_{1}+x_{2}>0 ,\\ \sin \left(x_{1}+x_{2}\right)+\cos \left(x_{1}+x_{2}\right), & \text { if } x_{1}+x_{2}\leq 0.\end{cases}
\end{array}\right.
\end{equation*}
The coefficient $\beta$ is
\begin{equation*}
\beta=\left\{\begin{array}{l}
\beta^{-}=\left(x_1^{2}-x_2^{2}+3\right) / 7,\boldsymbol{x}\in \Omega^-,\\
\beta^{+}=8,\boldsymbol{x}\in \Omega^+.
\end{array}\right.
\end{equation*}
The exact interface is the zero level set of the following level set function,
\begin{equation*}
\varphi(x)= \begin{cases}x_{2}-2 x_{1}, & \text { if } x_{1}+x_{2}>0, \\ x_{2}+x_{1} / 2, & \text { if } x_{1}+x_{2} \leq 0.\end{cases}
\end{equation*}
\begin{table}{}
 	\centering
 	\caption{$L^2$ errors and convergence orders of the DNN-FD method for Example \ref{example6}.}
 	\label{table10}
 	\begin{tabular}{llllllll}
 		& $\left\|u_{h}-u\right\|_{L^2(\Omega_1)}$ &   & $\left\|u_{h}-u\right\|_{L^2(\Omega_2)}$&  & $\left\|u_{h}-u\right\|_{L^2(\Omega)}$ &  & \\
 		N & Error & Order & Error & Order & Error & Order\\
 		\hline 20 $\times$ 20 & $2.48 \mathrm{E}-02$ & $-$ & $2.55 \mathrm{E}-02$ & $-$ & $1.99 \mathrm{E}-02$ & $-$ \\
 		40 $\times$ 40 & $5.89 \mathrm{E}-03$ & $2.0780$ & $6.13 \mathrm{E}-03$ & $2.0584$ & $5.52 \mathrm{E}-03$ & $1.8529$ \\
 		80 $\times$ 80 & $1.50 \mathrm{E}-03$ & $1.9669$ & $1.56 \mathrm{E}-03$ & $1.9718$ & $1.41 \mathrm{E}-03$ & $1.9680$ \\
 		160 $\times$ 160 & $3.85 \mathrm{E}-04$ & $1.9667$ & $3.99 \mathrm{E}-04$ & $1.9710$ & $3.53 \mathrm{E}-04$ & $1.9982$ \\
 		320 $\times$ 320 & $9.38 \mathrm{E}-05$ & $2.0388$ & $9.742 \mathrm{E}-05$ & $2.0348$ & $9.30 \mathrm{E}-05$ & $1.9246$
 	\end{tabular}
\end{table}
\begin{table}{}
	\centering
	\caption{$L^\infty$ errors and convergence orders of the DNN-FD method for Example \ref{example6}.}
	\label{table13}
	\begin{tabular}{lllllll}
		& $\left\|u_{h}-u\right\|_{L^\infty(\Omega)}$ & & IFVE\cite{wang2021new}\\
		N  & Error & Order & Order \\
		\hline 20 $\times$ 20 & $1.51 \mathrm{E}-02$ & $-$  & $-$\\
		40 $\times$ 40 &  $3.64 \mathrm{E}-03$ & $2.0544$ & $1.3132$\\
		80 $\times$ 80 &  $1.09 \mathrm{E}-03$ & $1.7943$ & $1.0505$\\
		160 $\times$ 160 &  $3.06 \mathrm{E}-04$ & $1.7793$ & $1.0106$\\
		320 $\times$ 320 & $8.96 \mathrm{E}-05$ & $1.7716$& $1.0139$
	\end{tabular}
\end{table}
\begin{figure}
	\centering    
	\subfigure[DNN-FD Solution]
	{
		\begin{minipage}{5cm}
			\centering
			\includegraphics[width=5cm,clip]{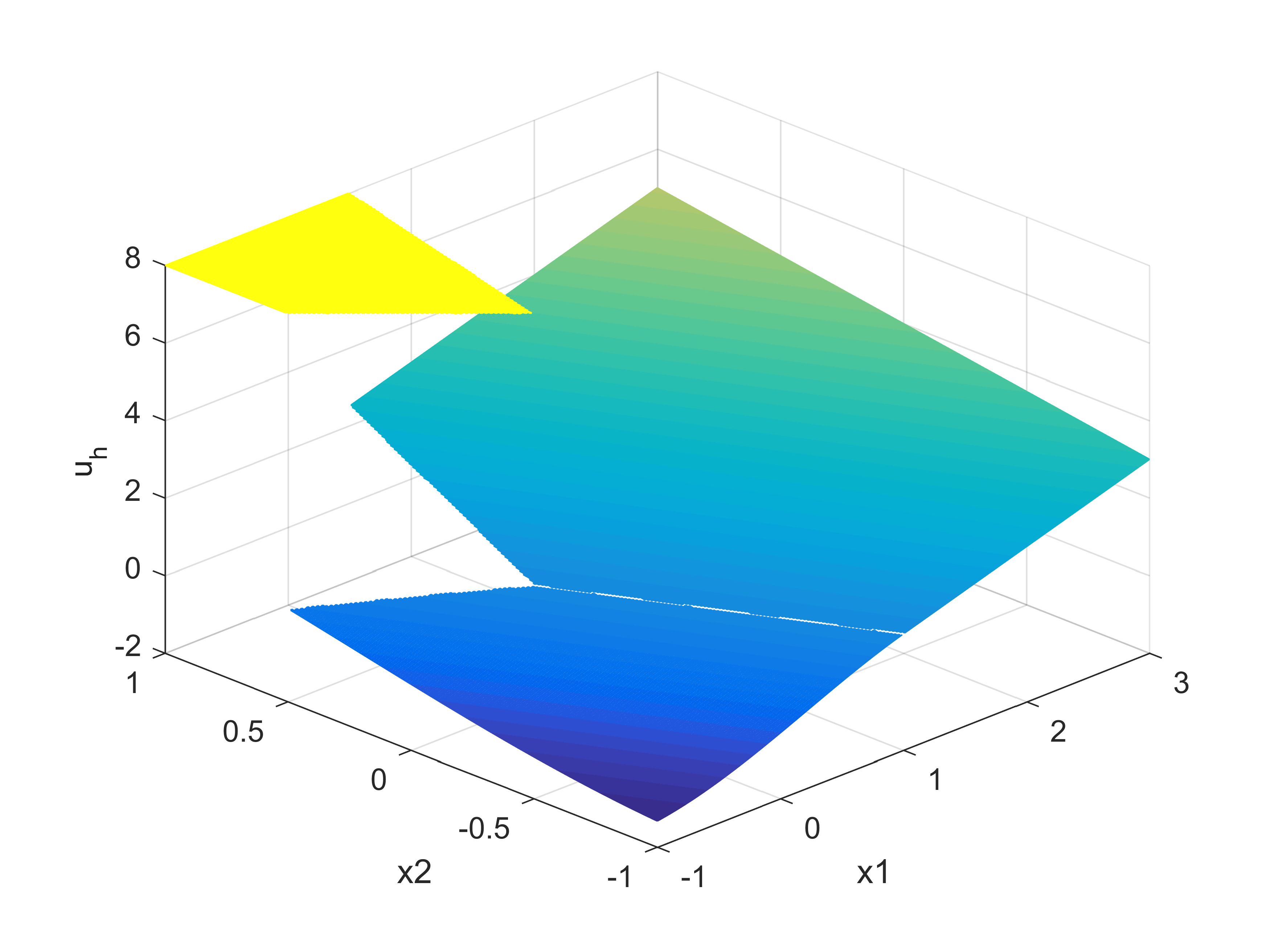}
			
		\end{minipage}
	}
	\subfigure[Exact Solution]
	{
		\begin{minipage}{5cm}
			\centering
			\includegraphics[width=5cm,clip]{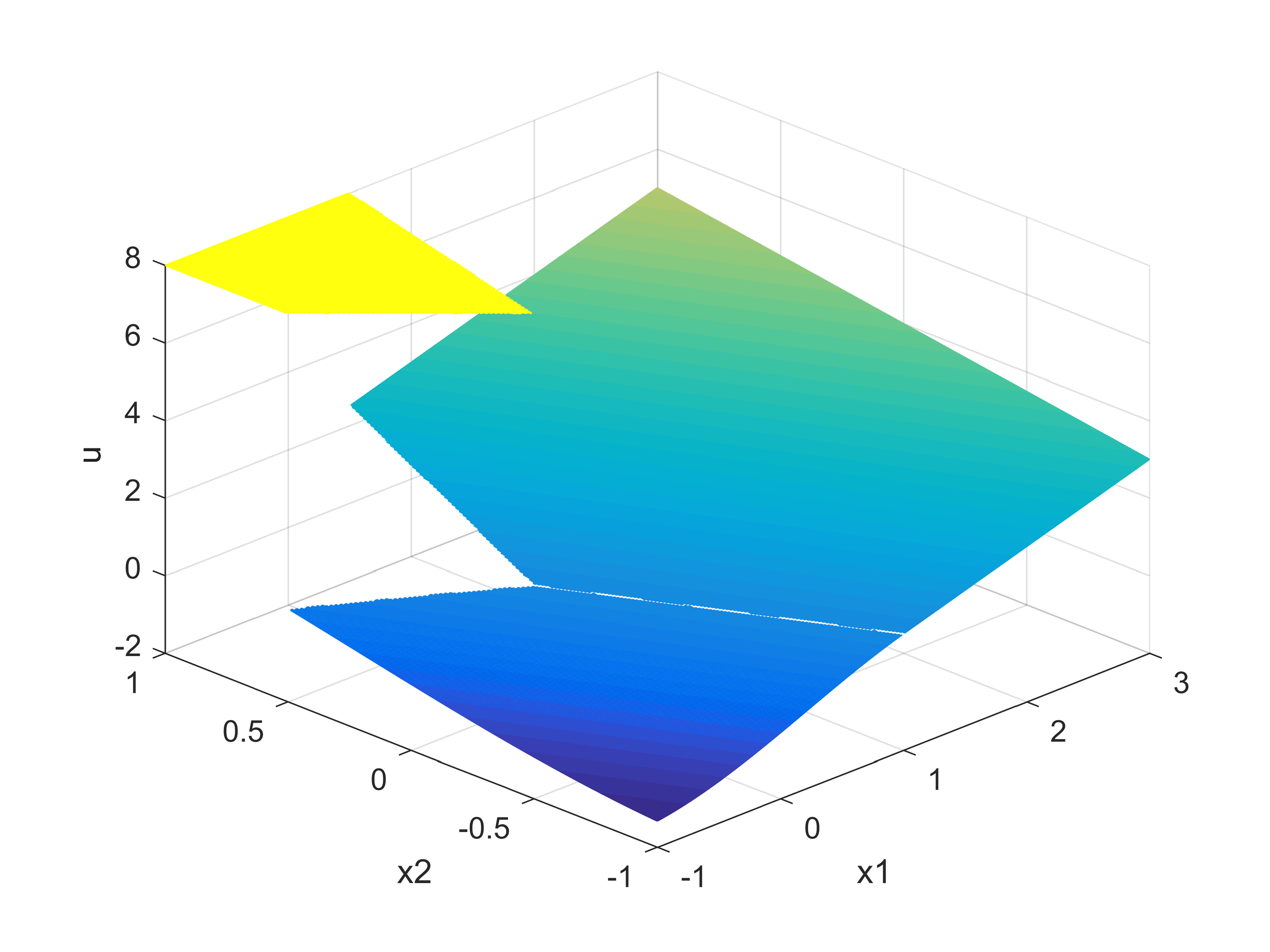}
			
		\end{minipage}
	}
	\vspace{-3mm}
	\caption{Comparison between exact and DNN-FD solutions for Example \ref{example6} when N=320.}
	\label{figure14}
\end{figure}
For nonsmooth interface problems, the method used in this paper can also be applied the numerical results of the current method are given in Table \ref{table10}. In Table \ref{table10}, we present a grid refinement analysis that successfully achieves the second order. In other words, the proposed method is not sensitive to the grid for the solution and interface. In Table \ref{table13}, We also calculated the logarithmic ratios of $L^\infty$ errors. Although the scheme is the second order one and costs too much expensive works on the interface, it is so hard to get satisfactory results in \cite{wang2021new} because of nonsmooth property of the interface. And the solution $u$ has a singularity at $(0, 0)$ with blow-up derivatives. Our method has approximately the second-order convergence, the numercial results are much better than ones of IFVE method. Fig.\ref{figure14} shows the comparison between the exact solution and the numerical solution when N=320.


\subsection{2D nondegenerate five-pointed star interface}
\label{example7}
\textbf{Example 4.8.} In this example, we consider the five-pointed star interface problem. The exact solution is\cite{hou2010numerical}
\begin{equation*}
u(\boldsymbol{x})=\left\{\begin{array}{l}
u^-(\boldsymbol{x})=8,\boldsymbol{x}\in \Omega^-, \\
u^+(\boldsymbol{x})=x_1^{2}+x_2^{2}+\sin (x_1+x_2),\boldsymbol{x}\in \Omega^+.
\end{array}\right.
\end{equation*}
The coefficient $\beta$ is
\begin{equation*}
\beta=\left\{\begin{array}{l}
\beta^{-}=1,\boldsymbol{x}\in \Omega^-,\\
\beta^{+}=2+\sin (x_1+x_2),\boldsymbol{x}\in \Omega^+.
\end{array}\right.
\end{equation*}
The exact interface is the zero level set of the following level set function,
\begin{equation*}
\phi(r, \theta)=\left\{\begin{array}{l}
\frac{R \sin \left(\theta_{t} / 2\right)}{\sin \left(\theta_{t} / 2+\theta-\theta_{r}-2 \pi(i-1) / 5\right)}-r,
\theta_{r}+\frac{\pi(2 i-2)}{5} \leqslant \theta<\theta_{r}+\frac{\pi(2 i-1)} {5},\\
\frac{R \sin \left(\theta_{t} / 2\right)}{\sin \left(\theta_{t} / 2-\theta+\theta_{r}-2 \pi(i-1) / 5\right)}-r,
\theta_{r}+\frac{\pi(2 i-3)}{5} \leqslant \theta<\theta_{r}+\frac{\pi(2 i-2)}{5} .
\end{array}\right.
\end{equation*}
with $\theta_{t}=\pi / 5, \theta_{r}=\pi / 7, R=6 / 7 \text { and } i=1,2,3,4,5 .$
\begin{table}{}
	\centering
	\caption{$L^2$ errors and convergence orders of the DNN-FD method for Example \ref{example7}.}
	\label{table11}
	\begin{tabular}{lllllll}
		& $\left\|u_{h}-u\right\|_{L^2(\Omega_1)}$ &   & $\left\|u_{h}-u\right\|_{L^2(\Omega_2)}$&  & $\left\|u_{h}-u\right\|_{L^2(\Omega)}$\\
		N & Error & Order & Error & Order & Error & Order \\
		\hline 20 $\times$ 20 & $2.57 \mathrm{E}-02$ & $-$ & $1.39 \mathrm{E}-02$ & $-$ & $2.09 \mathrm{E}-02$ & $-$ \\
		40 $\times$ 40 & $6.13 \mathrm{E}-03$ & $2.0584$ & $3.54 \mathrm{E}-03$ & $1.9708$ & $5.84 \mathrm{E}-03$ & $1.8436$ \\
		80 $\times$ 80 & $1.56 \mathrm{E}-03$ & $1.9718$ & $8.57 \mathrm{E}-04$ & $2.0487$ & $1.48 \mathrm{E}-03$ & $1.9716$ \\
		160 $\times$ 160 & $3.99 \mathrm{E}-04$ & $1.9710$ & $2.17 \mathrm{E}-04$ & $1.9818$ & $3.74 \mathrm{E}-04$ & $1.9898$ \\
		320 $\times$ 320 & $9.74 \mathrm{E}-05$ & $2.0348$ & $5.36 \mathrm{E}-05$ & $2.0160$ & $9.79 \mathrm{E}-05$ & $1.9364$
	\end{tabular}
\end{table}
\begin{figure}
	\centering    
	\subfigure[DNN-FD Solution]
	{
		\begin{minipage}{5cm}
			\centering
			\includegraphics[width=5cm,clip]{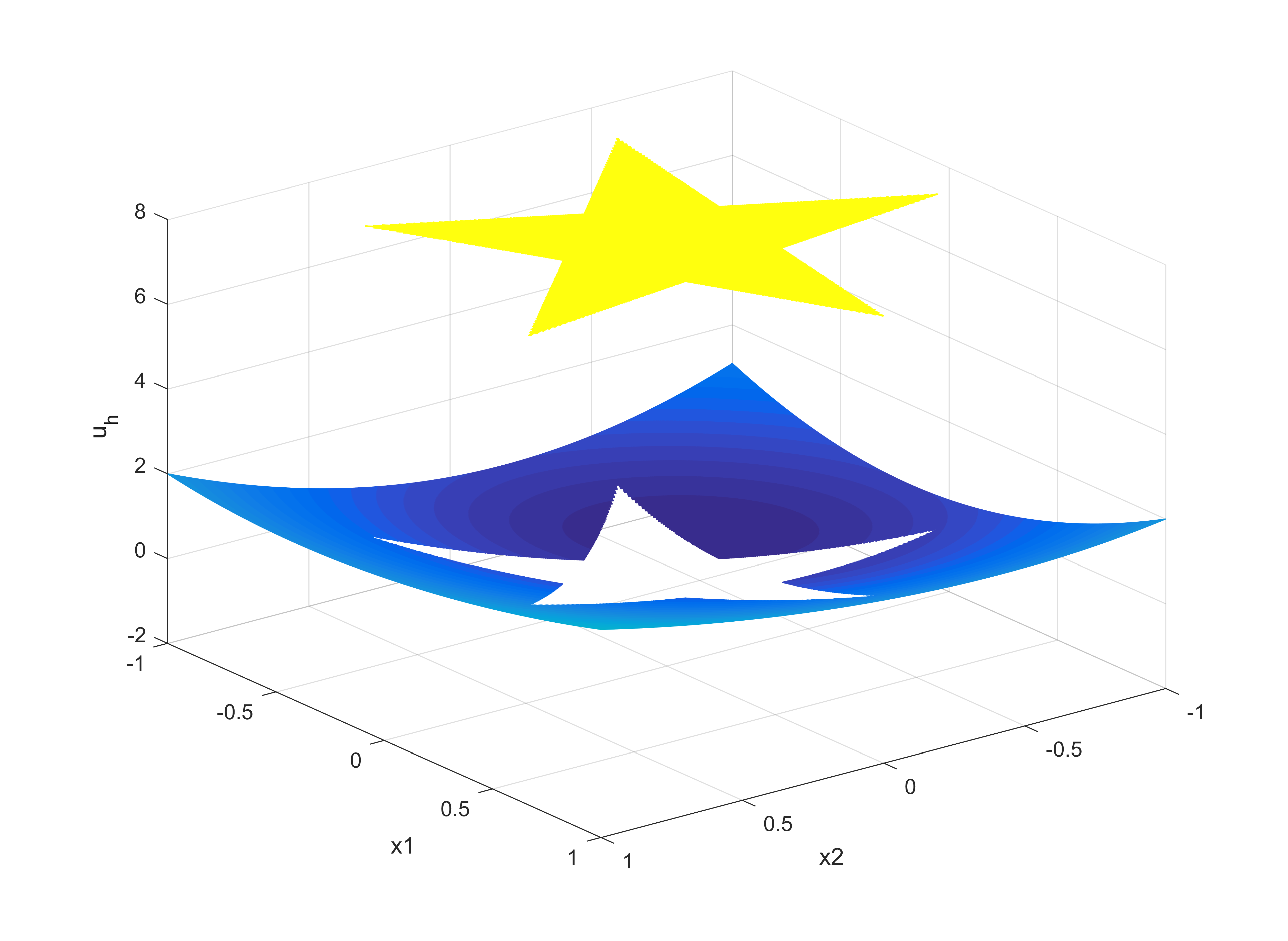}
			
		\end{minipage}
	}
	\subfigure[Exact Solution]
	{
		\begin{minipage}{5cm}
			\centering
			\includegraphics[width=5cm,clip]{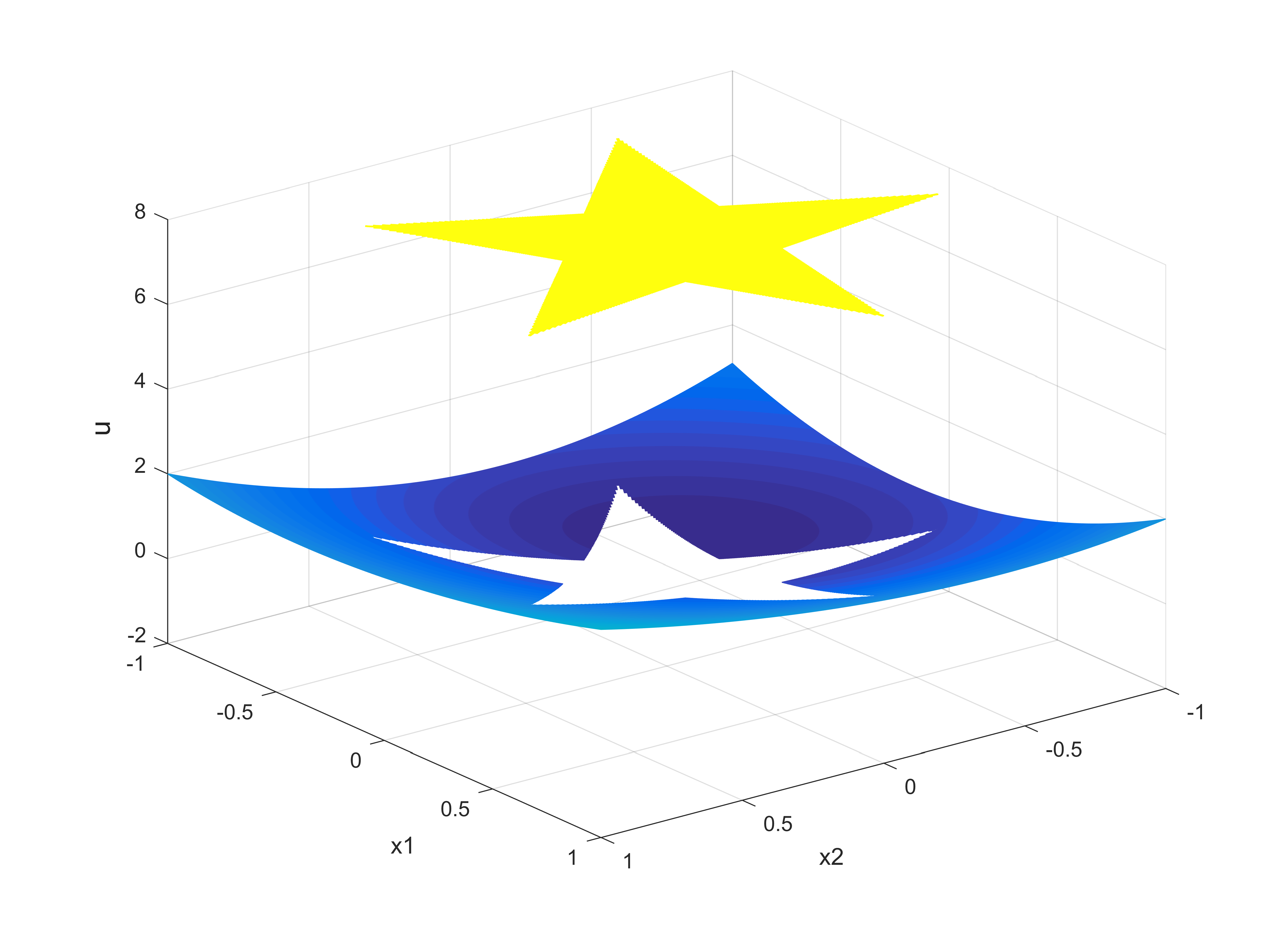}
			
		\end{minipage}
	}
	\vspace{-3mm}
	\caption{Comparison between exact and DNN-FD solutions for Example \ref{example7} when N=320.}
	\label{figure15}
\end{figure}
This example presents a more difficult challenge, that is, considering that the interface consists of several sharp-edged nonsmooth interfaces. Our method can also be applied after special processing for more complex nonsmooth interfaces, such as the five-pointed star interface. The numerical results of the current method for in Table \ref{table11}. It can be seen that even if the non-smoothness of the interface changes, our method can always maintain the second-order accuracy. The exact solution and the numerical solution are compared in Fig.\ref{figure15} when N=320.


\subsection{2D degenerate five-pointed star interface}
\label{example8}
\textbf{Example 4.9.} In this example, we consider the degenerate five-pointed star interface problem. The exact solution is\cite{hou2010numerical}
\begin{equation*}
u(\boldsymbol{x})=\left\{\begin{array}{l}
u^-(\boldsymbol{x})=6+\sin (2 \pi x_1) \sin (2 \pi x_2),\boldsymbol{x}\in \Omega^-,\\
u^+(\boldsymbol{x})=x_1^{2}+x_2^{2}+\sin (x_1+x_2),\boldsymbol{x}\in \Omega^+.
\end{array}\right.
\end{equation*}
The coefficient $\beta$ is
\begin{equation*}
\beta=\left\{\begin{array}{l}
\beta^{-}=(x_1-\frac{6}{7})^2+(x_2-\frac{6}{7})^2,\boldsymbol{x}\in \Omega^-,\\
\beta^{+}=(x_1-\frac{6sin(\pi/10)}{7sin(\pi/3)})^2+(x_2-\frac{6sin(\pi/10)}{7sin(\pi/3)})^2,\boldsymbol{x}\in \Omega^+.
\end{array}\right.
\end{equation*}

The exact interface is the same as in the previous example.
\begin{table}{}
	\centering
	\caption{$L^2$ errors and convergence orders of the DNN-FD method for Example \ref{example8}.}
	\label{table12}
	\begin{tabular}{lllllll}
		& $\left\|u_{h}-u\right\|_{L^2(\Omega_1)}$ &   & $\left\|u_{h}-u\right\|_{L^2(\Omega_2)}$&  & $\left\|u_{h}-u\right\|_{L^2(\Omega)}$\\
		N & Error & Order & Error & Order & Error & Order \\
		\hline 20 $\times$ 20 & $1.30 \mathrm{E}-02$ & $-$ & $1.50 \mathrm{E}-02$ & $-$ & $1.88 \mathrm{E}-02$ & $-$ \\
		40 $\times$ 40 & $3.54 \mathrm{E}-03$ & $1.8772$ & $3.63 \mathrm{E}-03$ & $2.0499$ & $5.10 \mathrm{E}-03$ & $1.8842$ \\
		80 $\times$ 80 & $7.97 \mathrm{E}-04$ & $2.1508$ & $1.42 \mathrm{E}-03$ & $1.3528$ & $1.24 \mathrm{E}-03$ & $2.0385$ \\
		160 $\times$160 & $2.49 \mathrm{E}-04$ & $1.6769$ & $3.64 \mathrm{E}-04$ & $1.9644$ & $3.02 \mathrm{E}-04$ & $2.0391$ \\
		320 $\times$ 320 & $6.11 \mathrm{E}-05$ & $2.0268$ & $9.25 \mathrm{E}-05$ & $1.9795$ & $8.11 \mathrm{E}-05$ & $1.8973$
	\end{tabular}
\end{table}
\begin{figure}
	\centering    
	\subfigure[DNN-FD Solution]
	{
		\begin{minipage}{5cm}
			\centering
			\includegraphics[width=5cm,clip]{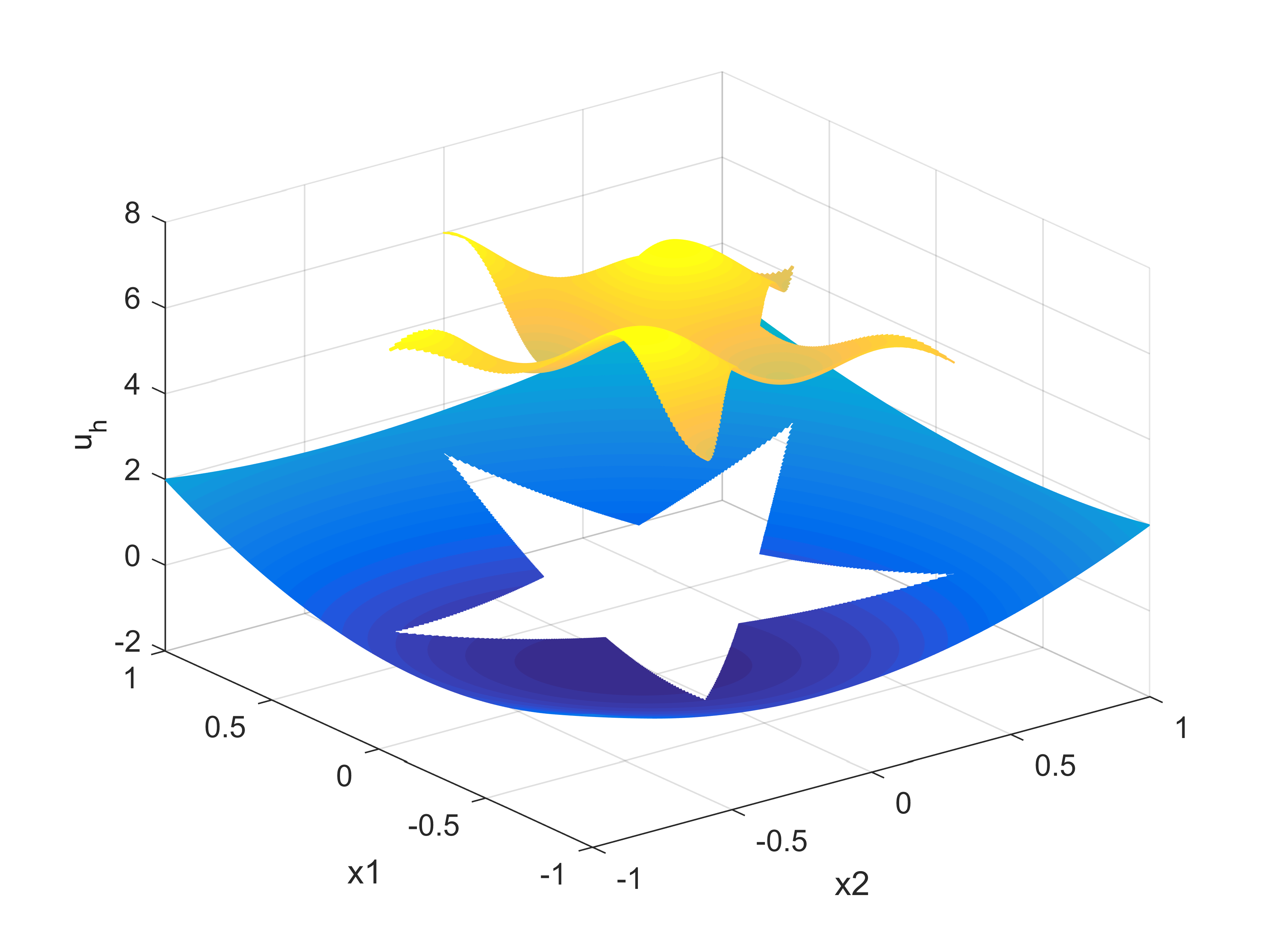}
			
		\end{minipage}
	}
	\subfigure[Exact Solution]
	{
		\begin{minipage}{5cm}
			\centering
			\includegraphics[width=5cm,clip]{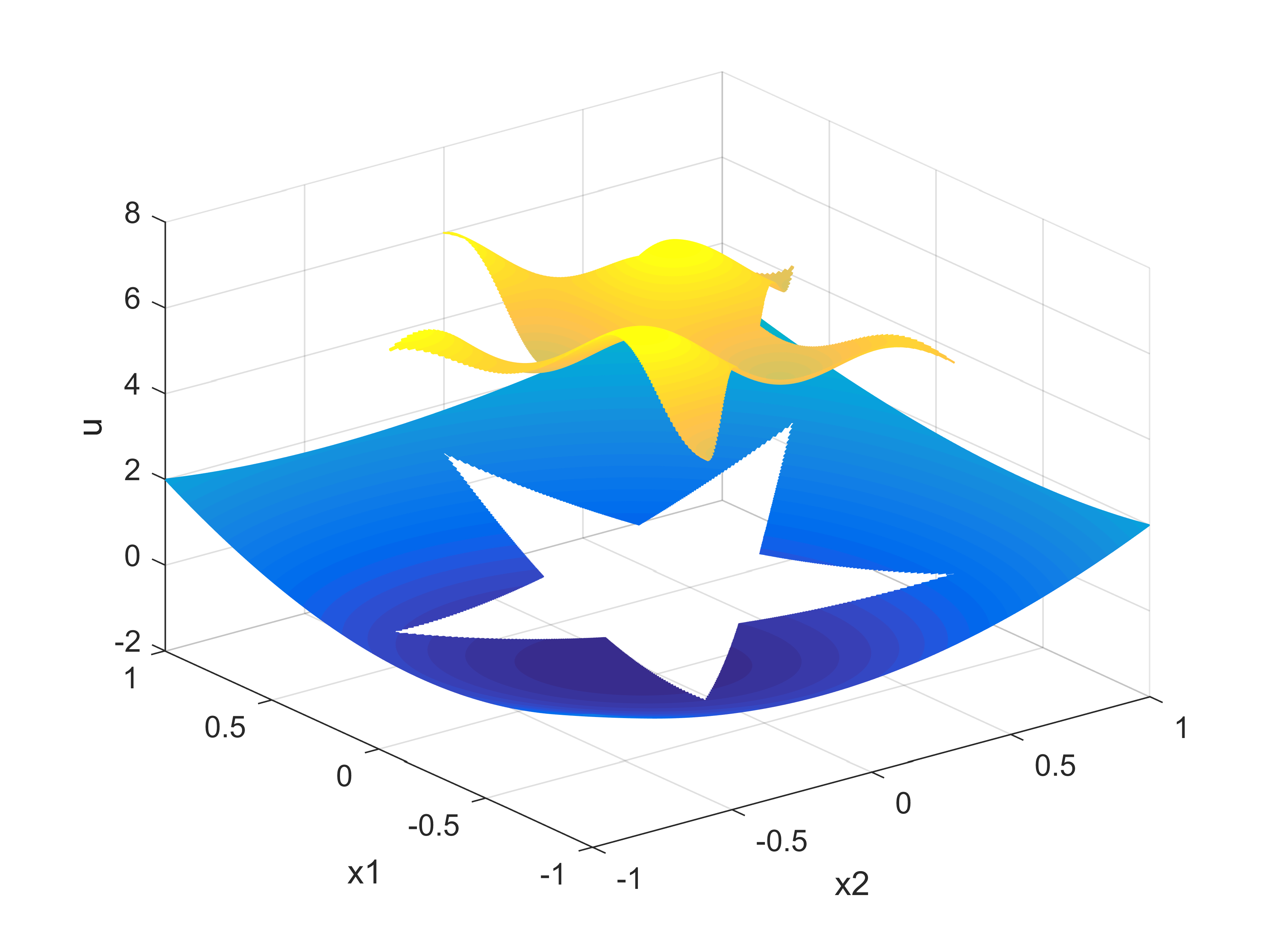}
			
		\end{minipage}
	}
	\vspace{-3mm}
	\caption{Comparison between exact and DNN-FD solutions for Example \ref{example8} when N=320.}
	\label{figure16}
\end{figure}
In the last example, we reconstruct the examples from the original literature\cite{hou2010numerical}. We will challenge one which is combining degenerate and nonsmooth interface problems, where the degenerate points are respectively the two angles of the five-pointed star on the positive and negative domains. Furthermore, because the solution of the problem is nonlinear, the difficulty of this example increases once again. The choice of the activation function has also changed, and the selected nonlinear activation function offers a good approximation to the solution of the problem. The numerical results of the current method are shown in Table \ref{table12}. The experimental results have the second-order accuracy in the  $L^2$ norm. Fig.\ref{figure16} shows the comparison between the exact solution and the numerical solution when N=320.

\subsection{2D degenerate interface with large jump conditions}
\label{example10}
\textbf{Example 4.10.} This example is based on the addition of a large jump ratio to Example \ref{example8}. The boundary condition and the source function are chosen so that the exact solution is\cite{zhao2021semi}
\begin{equation*}
u(x)=\left\{
\begin{array}{l}
7 x_1^{2}+7 x_2^{2}+6,x \in \Omega^{-}, \\
x_1^{2}+x_2^{2}+\sin (x_1+x_2), x \in \Omega^{+}.
\end{array}\right.
\end{equation*}
The coefficient $\beta$ is
\begin{equation*}
\beta=\left\{\begin{array}{l}
\beta^{-}=\tau^{-}((x_1-\frac{6}{7})^2+(x_2-\frac{6}{7})^2),\boldsymbol{x}\in \Omega^-,\\
\beta^{+}=\tau^{+}((x_1-\frac{6sin(\pi/10)}{7sin(\pi/3)})^2+(x_2-\frac{6sin(\pi/10)}{7sin(\pi/3)})^2),\boldsymbol{x}\in \Omega^+.
\end{array}\right.
\end{equation*}
The exact interface is the zero level set of the following level set function,
\begin{equation*}
\phi(r, \theta)=\left\{\begin{array}{l}
\frac{R \sin \left(\theta_{t} / 2\right)}{\sin \left(\theta_{t} / 2+\theta-\theta_{r}-2 \pi(i-1) / 5\right)}-r,
\theta_{r}+\frac{\pi(2 i-2)}{5} \leqslant \theta<\theta_{r}+\frac{\pi(2 i-1)} {5},\\
\frac{R \sin \left(\theta_{t} / 2\right)}{\sin \left(\theta_{t} / 2-\theta+\theta_{r}-2 \pi(i-1) / 5\right)}-r,
\theta_{r}+\frac{\pi(2 i-3)}{5} \leqslant \theta<\theta_{r}+\frac{\pi(2 i-2)}{5} .
\end{array}\right.
\end{equation*}
with $\theta_{t}=\pi / 5, \theta_{r}=\pi / 7, R=6 / 7 \text { and } i=1,2,3,4,5 .$
\begin{figure}[htp]
	\centering
	\includegraphics[width=10cm]{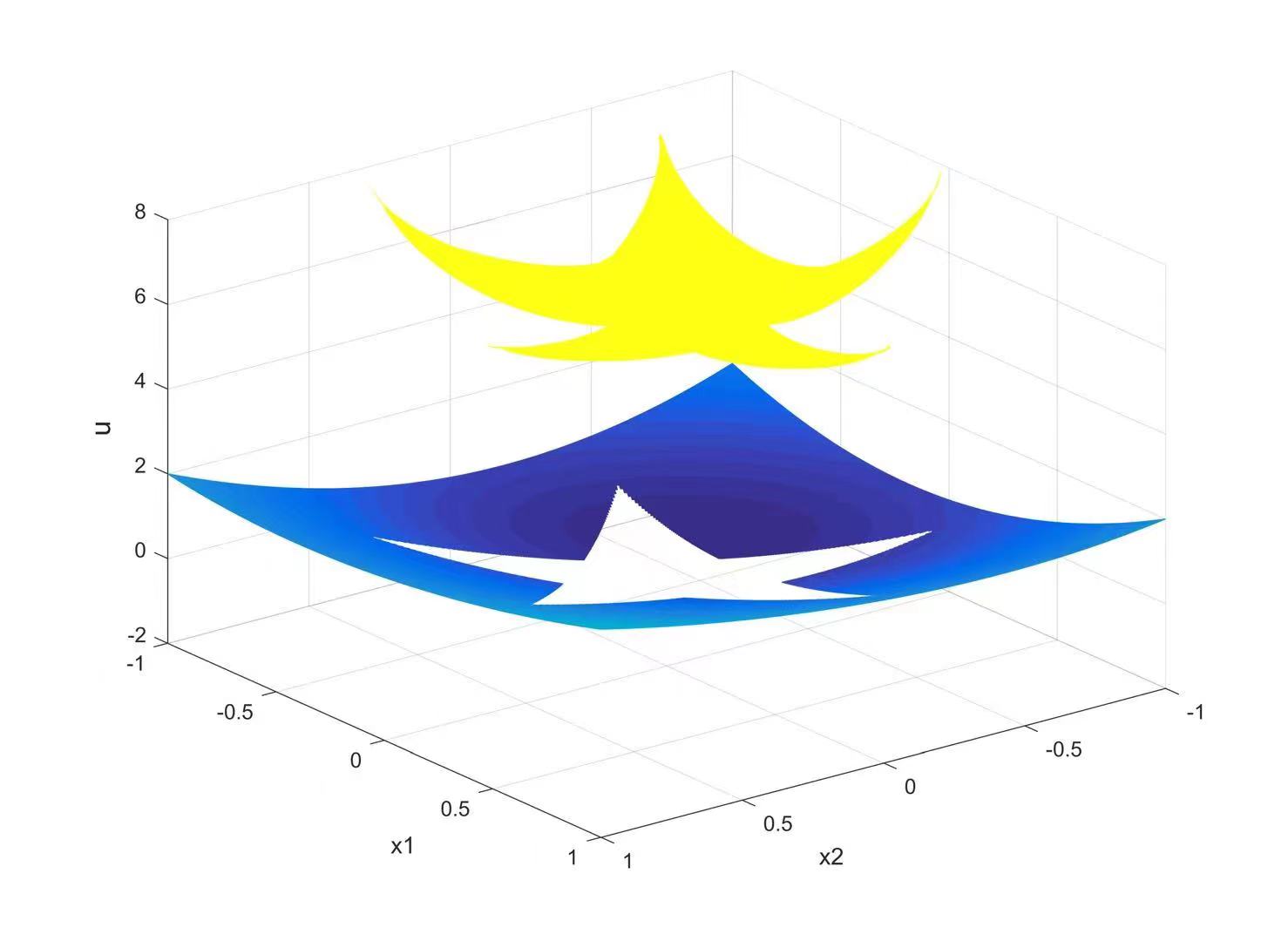}\\
	\caption{The DNN-FD solution for Example \ref{example10} when N=320 ($\tau^{-}/\tau^{+}=1:10^{10}$).}
	\label{figure22}
\end{figure}
\begin{table}{}
	\centering
	\caption{$L^2$ errors and convergence orders of the DNN-FD method for Example \ref{example10} ($\tau^{-}/\tau^{+}=1:10^{10}$).}
	\label{table18}
	\begin{tabular}{lllllll}
		& $\left\|u_{h}-u\right\|_{L^2(\Omega_1)}$ &   & $\left\|u_{h}-u\right\|_{L^2(\Omega_2)}$&  & $\left\|u_{h}-u\right\|_{L^2(\Omega)}$\\
		N & Error & Order & Error & Order & Error & Order \\
		\hline 20 $\times$ 20 & $1.30 \mathrm{E}-02$ & $-$ & $1.50 \mathrm{E}-02$ & $-$ & $1.88 \mathrm{E}-02$ & $-$ \\
		40 $\times$ 40 & $3.54 \mathrm{E}-03$ & $1.8772$ & $3.63 \mathrm{E}-03$ & $2.0499$ & $5.10 \mathrm{E}-03$ & $1.9842$ \\
		80 $\times$ 80 & $7.67 \mathrm{E}-04$ & $2.1508$ & $1.42 \mathrm{E}-03$ & $1.3528$ & $1.24 \mathrm{E}-03$ & $2.0325$ \\
		160 $\times$160 & $1.49 \mathrm{E}-04$ & $1.6769$ & $3.64 \mathrm{E}-04$ & $1.9644$ & $3.02 \mathrm{E}-04$ & $1.8991$ \\
		320 $\times$ 320 & $5.12 \mathrm{E}-05$ & $2.0268$ & $9.25 \mathrm{E}-05$ & $1.9795$ & $8.11 \mathrm{E}-05$ & $1.8273$
	\end{tabular}
\end{table}
\begin{table}{}
	\centering
	\caption{$L^2$ errors and convergence orders of the DNN-FD method for Example \ref{example10} ($\tau^{-}/\tau^{+}=10^{10}:1$).}
	\label{table19}
	\begin{tabular}{lllllll}
		& $\left\|u_{h}-u\right\|_{L^2(\Omega_1)}$ &   & $\left\|u_{h}-u\right\|_{L^2(\Omega_2)}$&  & $\left\|u_{h}-u\right\|_{L^2(\Omega)}$\\
		N & Error & Order & Error & Order & Error & Order \\
		\hline 20 $\times$ 20 & $1.29 \mathrm{E}-02$ & $-$ & $1.50 \mathrm{E}-02$ & $-$ & $1.88 \mathrm{E}-02$ & $-$ \\
		40 $\times$ 40 & $3.54 \mathrm{E}-03$ & $1.8772$ & $3.63 \mathrm{E}-03$ & $2.0499$ & $5.10 \mathrm{E}-03$ & $1.8842$ \\
		80 $\times$ 80 & $7.97 \mathrm{E}-04$ & $2.1508$ & $5.12 \mathrm{E}-03$ & $1.3528$ & $1.24 \mathrm{E}-03$ & $1.9385$ \\
		160 $\times$160 & $2.49 \mathrm{E}-04$ & $1.6769$ & $3.62 \mathrm{E}-04$ & $1.9644$ & $3.02 \mathrm{E}-04$ & $1.8591$ \\
		320 $\times$ 320 & $6.11 \mathrm{E}-05$ & $2.0268$ & $8.15 \mathrm{E}-05$ & $1.9795$ & $8.11 \mathrm{E}-05$ & $1.7903$
	\end{tabular}
\end{table}

Our method can also be applied in the five-pointed star interface with large jump ratios. The numerical results of the current method for in Table \ref{table18} and Table\ref{table19}. It can be seen that even if the non-smoothness of the interface changes, our method can always maintain the second-order accuracy. The numerical solution is shown in Fig.\ref{figure22} when N=320.

\subsection{2D interface problem with non analytical solution}
\label{example11}
\textbf{Example 4.11.} In this example, we consider the five-pointed star interface problem with non analytical solution which is constructed from Example \ref{example7}.
The coefficient $\beta$ is
\begin{equation*}
\beta=\left\{\begin{array}{l}
\beta^{-}=1,\boldsymbol{x}\in \Omega^-,\\
\beta^{+}=2+\sin (x_1+x_2),\boldsymbol{x}\in \Omega^+.
\end{array}\right.
\end{equation*}
The exact interface is the zero level set of the following level set function,
\begin{equation*}
\phi(r, \theta)=\left\{\begin{array}{l}
\frac{R \sin \left(\theta_{t} / 2\right)}{\sin \left(\theta_{t} / 2+\theta-\theta_{r}-2 \pi(i-1) / 5\right)}-r,
\theta_{r}+\frac{\pi(2 i-2)}{5} \leqslant \theta<\theta_{r}+\frac{\pi(2 i-1)} {5},\\
\frac{R \sin \left(\theta_{t} / 2\right)}{\sin \left(\theta_{t} / 2-\theta+\theta_{r}-2 \pi(i-1) / 5\right)}-r,
\theta_{r}+\frac{\pi(2 i-3)}{5} \leqslant \theta<\theta_{r}+\frac{\pi(2 i-2)}{5} .
\end{array}\right.
\end{equation*}
with $\theta_{t}=\pi / 5, \theta_{r}=\pi / 7, R=6 / 7 \text { and } i=1,2,3,4,5 .$
\begin{table}{}
	\centering
	\caption{ $L^2$ errors and convergence orders of $f$ and $f_h$  for Example \ref{example11}.}
	\label{table17}
	\begin{tabular}{lllllll}
		& $\left\|f_{h}-f\right\|_{L^2(\Omega_1)}$ &   & $\left\|f_{h}-f\right\|_{L^2(\Omega_2)}$&  & $\left\|f_{h}-f\right\|_{L^2(\Omega)}$\\
		N & Error & Order & Error & Order & Error & Order \\
		\hline 20 $\times$ 20 & $2.57 \mathrm{E}-02$ & $-$ & $1.39 \mathrm{E}-02$ & $-$ & $2.09 \mathrm{E}-02$ & $-$ \\
		40 $\times$ 40 & $6.13 \mathrm{E}-03$ & $1.5584$ & $3.54 \mathrm{E}-03$ & $1.4708$ & $5.84 \mathrm{E}-03$ & $1.4436$ \\
		80 $\times$ 80 & $1.56 \mathrm{E}-03$ & $1.3718$ & $8.57 \mathrm{E}-04$ & $1.5487$ & $1.48 \mathrm{E}-03$ & $1.4716$ \\
		160 $\times$ 160 & $3.99 \mathrm{E}-04$ & $1.4710$ & $2.17 \mathrm{E}-04$ & $1.4818$ & $3.74 \mathrm{E}-04$ & $1.5898$ \\
		320 $\times$ 320 & $9.74 \mathrm{E}-05$ & $1.4348$ & $5.36 \mathrm{E}-05$ & $1.5160$ & $9.79 \mathrm{E}-05$ & $1.5364$
	\end{tabular}
\end{table}
\begin{figure}[htp]
	\centering
	\includegraphics[width=10cm]{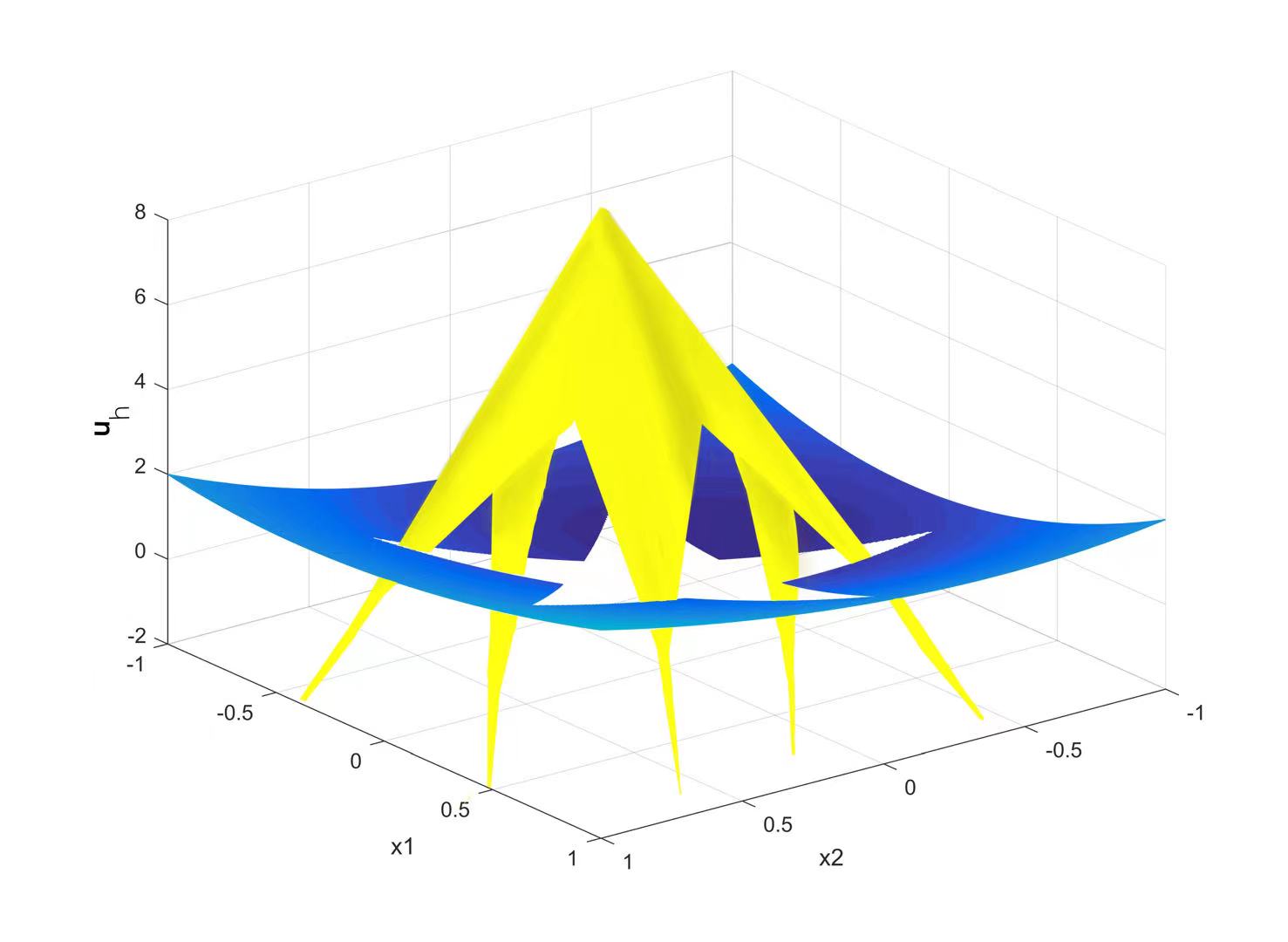}\\
	\caption{The DNN-FD solution for Example \ref{example11}. when N=320.}
	\label{figure21}
\end{figure}
We change the $f^-(\boldsymbol{x})=|\boldsymbol{x}-\boldsymbol{x}_0|(1+2 \log |\boldsymbol{x}-\boldsymbol{x}_0|)$, where $\phi(\boldsymbol{x}_0)=0$. This example presents a more difficult challenge, that is, considering that the interface consists of several sharp-edged nonsmooth interfaces and has non analytical solution. Our method can also be applied this example. The numerical results of the current method for in Table \ref{table17} where $f_h$ is the right term calculated by the numerical solution $u_h$. Due to the lack of the analytical solution to the equation, we define $L^2$ errors and convergence orders of the equation as the reference of stability during the operation. This value is stable around a constant, confirming the feasibility of the method. The numerical solution is shown in Fig.\ref{figure21} when N=320.

\subsection{2D Linear elasticity interface problem}
\label{example12}
\textbf{Example 4.12.} Finally, we will consider the example with physical significance that is a linear elasticity PDE with a discontinuous stress tensor as follows,
\begin{equation*}
-\nabla \cdot \mathbb{T}=f(\boldsymbol{x},u), \text { in } \Omega^- \cup \Omega^+,
\end{equation*}
\begin{equation*}
[u]=w, \text { on } \Gamma,
\end{equation*}
\begin{equation*}
[\mathbb{T}  \cdot \boldsymbol{n}]=v, \text { on } \Gamma,
\end{equation*}
\begin{equation*}
u=g, \text { on } \partial \Omega.
\end{equation*}
One application of the linear elasticity problem is to model the shape and location of fibroblast cells under stress. Let $\mathbf{u}=\left(u_1, u_2\right)^T$ denote the displacement field. Then, the strain tensor is
\begin{equation*}
\sigma=\frac{1}{2}\left(\nabla \mathbf{u}+(\nabla \mathbf{u})^T\right).
\end{equation*}
then the elasticity tensor $\mathbb{T}$ is a linear transformation on the tensors. In the isotropic case, we have
\begin{equation*}
\mathbb{T}\sigma=\lambda \operatorname{Tr}(\sigma) \mathbf{1}+2\mu\left(\sigma+\sigma^T\right).
\end{equation*}
where $\lambda$ and $\mu$ are lamé constants, $\operatorname{Tr(.)}$ is the trace operator, and $\mathbf{1}$ is the identity matrix. In this case, the above parameters satisfies the following relationships
\begin{equation*}
\mu=\frac{E}{2(1+\nu)}, \quad \lambda=\frac{E \nu}{(1+\nu)(1-2 \nu)}.
\end{equation*}
where E is Young modulus and $\mu,\nu$ are Poisson’s ratio. The
interface is  defined in the polar coordinate
\begin{equation*}
r=0.5+\frac{\sin 5 \theta}{7}.
\end{equation*}
We set the computational domain $\Omega =[-1,1]\times[-1,1].$ The Dirichlet boundary condition and homogeneous jump conditions are determined in this example. Then we choose two groups of the Poisson’s ratio and the shear modulus as follows\cite{wang2015matched}
\begin{equation}
\label{e7}
\nu= \begin{cases}\nu^{-}=0.24, & \text { in } \Omega^{-}; \\ \nu^{+}=0.20, & \text { in } \Omega^{+}.\end{cases},
\mu= \begin{cases}\mu^{-}=2000000, & \text { in } \Omega^{-}; \\ \mu^{+}=1500000, & \text { in } \Omega^{+}.\end{cases}
\end{equation}
and
\begin{equation}
\label{e8}
\nu= \begin{cases}\nu^{-}=0.24, & \text { in } \Omega^{-}; \\ \nu^{+}=0.00024, & \text { in } \Omega^{+} .\end{cases},
\mu= \begin{cases}\mu^{-}=2000000, & \text { in } \Omega^{-}; \\ \mu^{+}=1500000, & \text { in } \Omega^{+}.\end{cases}
\end{equation}

\begin{figure}[htp]
	\centering
	\includegraphics[width=10cm]{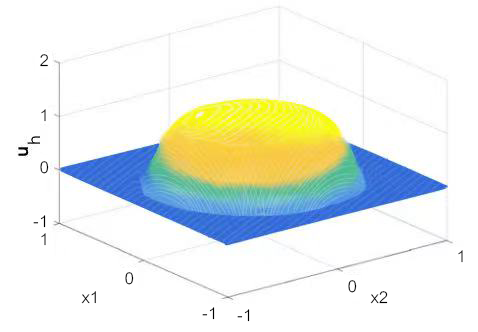}\\
	\caption{The DNN-FD solution for Example \ref{example12} with (\ref{e7}) when N=320.}
	\label{figure19}
\end{figure}

\begin{table}{}
	\centering
	\caption{$L^2$ errors and convergence orders with (\ref{e7}) for Example \ref{example12}.}
	\label{table15}
	\begin{tabular}{lllllll}
		& $\left\|u_{h}-u\right\|_{L^2(\Omega_1)}$ &   & $\left\|u_{h}-u\right\|_{L^2(\Omega_2)}$&  & $\left\|u_{h}-u\right\|_{L^2(\Omega)}$\\
		N & Error & Order & Error & Order & Error & Order \\
		\hline 20 $\times$ 20 & $5.99 \mathrm{E}-03$ & $-$ & $1.60 \mathrm{E}-03$ & $-$ & $3.39 \mathrm{E}-03$ & $-$ \\
		40 $\times$ 40 & $1.45 \mathrm{E}-04$ & $1.9267$ & $4.77 \mathrm{E}-04$ & $1.7577$ & $8.92 \mathrm{E}-04$ & $1.9234$ \\
		80 $\times$ 80 & $3.66 \mathrm{E}-04$ & $1.9746$ & $1.42 \mathrm{E}-05$ & $1.7503$ & $1.99 \mathrm{E}-04$ & $2.1647$ \\
		160 $\times$ 160 & $1.04 \mathrm{E}-04$ & $1.9919$ & $2.83 \mathrm{E}-05$ & $2.3223$ & $4.79 \mathrm{E}-05$ & $2.0076$ \\
		320 $\times$ 320 & $2.50 \mathrm{E}-05$ & $2.0653$ & $7.61 \mathrm{E}-06$ & $1.8912$ & $1.34\mathrm{E}-05$ & $1.8996$
	\end{tabular}
\end{table}

\begin{table}{}
	\centering
	\caption{$L^2$ errors and convergence orders with (\ref{e8}) for Example \ref{example12}.}
	\label{table16}
	\begin{tabular}{lllllll}
		& $\left\|u_{h}-u\right\|_{L^2(\Omega_1)}$ &   & $\left\|u_{h}-u\right\|_{L^2(\Omega_2)}$&  & $\left\|u_{h}-u\right\|_{L^2(\Omega)}$\\
		N & Error & Order & Error & Order & Error & Order \\
		\hline 20 $\times$ 20 & $4.08 \mathrm{E}-03$ & $-$ & $1.82 \mathrm{E}-03$ & $-$ & $1.38 \mathrm{E}-03$ & $-$ \\
		40 $\times$ 40 & $1.01 \mathrm{E}-04$ & $2.0108$ & $3.25 \mathrm{E}-04$ & $2.4939$ & $3.20 \mathrm{E}-04$ & $2.1162$ \\
		80 $\times$ 80 & $2.23 \mathrm{E}-04$ & $2.1886$ & $8.59 \mathrm{E}-05$ & $2.2130$ & $6.35 \mathrm{E}-05$ & $2.3367$ \\
		160 $\times$ 160 & $5.17 \mathrm{E}-05$ & $2.1108$ & $2.36 \mathrm{E}-05$ & $1.9278$ & $1.50 \mathrm{E}-05$ & $2.0828$ \\
		320 $\times$ 320 & $1.40 \mathrm{E}-06$ & $1.8868$ & $5.26 \mathrm{E}-06$ & $2.1720$ & $3.53 \mathrm{E}-06$ & $2.0864$
	\end{tabular}
\end{table}

The network used 6 intermediate layers. The width of each layer is 20 and the learning rate $\eta$ is $5 \times 10^{-4}$.In Fig.\ref{figure19}, we plot the profiles of the DNN-FD solution, which are the displacements in $x_1$ and $x_2$ coordinates, respectively. The corresponding numerical results are shown in Table \ref{table15} and Table \ref{table16}. We find that the DNN-FD solutions have the second-order accuracy in the $L^2$ norm.


\section{Conclusions.}
\label{s5}
Numerical methods for solving nolinear degenerate interface problems is one of fundamental iusses in scientific computations, it is challenge to design effective and robust fully decoupled numerical method for such degenerate interface problems. In this paper, fully decoupled finite difference method based on deep neural network for solving degenerate interface problems including 1D and 2D cases is proposed. It is shown that we can adopt uniform grids to solve degenerate PDE with interface. There are no unknown augmented parameters in the discrete schemes, and no more extra conditions and works to be required for designing numerical approximation algorithms. In fact, some augmented variables is obtained by adopting DNN technique, the degenerate interface problem is completely decoupled two independent
to the case of other degenerate or singular problems. The accuracy of the proposed fully decoupled algorithms has been demonstrated by solving various examples including degenerate and nondegenerate cases. In particular, the fully decoupled properties of the algorithm make the method capable of easy handling the jump ratio from the case of semi-decoupling $\mathbf{BIG}$ jump (such as $10^{7}:1$ or $1:10^{7}$) to the case of fully decoupled $\mathbf{VERY}$ $\mathbf{BIG}$ jump (such as $10^{12}:1$ or $1:10^{12}$) conditions. An interesting typical sharp edge example with degenerate five-pointed star interface shows that our approach works very well for those very hard problems. Numerical examples confirm the effectiveness of the fully decoupled algorithms for solving degenerate interface problems.


\section*{Acknowledgments.}
This work is partially supported by the National Natural Science Foundation of China(grants No. 11971241).

~\\

\renewcommand\refname{References}
\bibliographystyle{abbrv}

\end{document}